\documentclass[final]{amsart}

\usepackage{pdfsync}
\usepackage{amsmath,amssymb}
\usepackage{enumerate}
\usepackage{xspace}
\usepackage{boxedminipage}
\usepackage{algorithmic}
\usepackage{listings}
\usepackage{tabularx}
\usepackage{subfigure}
\usepackage{tristan}

\numberwithin{equation}{section}
\usepackage{showkeys}
\setboolean{showtodo}{true}
\setboolean{shownotes}{true}
\setboolean{showchanges}{true}
\setboolean{usemathrsfs}{false}
\author{Jan Giesselmann} 
\address{Jan Giesselmann\newline
 Institute of Applied Analysis and Numerical Simulation\newline
University of Stuttgart\newline
Pfaffenwaldring 57\newline
D-70563 Stuttgart\newline
  Germany} 
\curraddr{}
\email{\linkedemail{jan.giesselmann@mathematik.uni-stuttgart.de}} 
\thanks{The authors were supported by the  the FP7-REGPOT project ``ACMAC--Archimedes Center for Modeling, Analysis and Computations" of the University of Crete (FP7-REGPOT-2009-1-245749).}

\author{Charalambos Makridakis}
\address{ %
  Charalambos Makridakis\newline
  Department of Applied Mathematics\newline %
  University of Crete\newline %
  GR-71409 Heraklion, Greece\newline
  and
  \newline
  Institute for Applied and Computational Mathematics\newline %
  Foundation for Research and Technology-Hellas\newline %
  Vasilika Vouton P.O.Box 1527\newline %
  GR-71110 Heraklion, Greece}
\curraddr{}
\email{\linkedemail{makr@tem.uoc.gr}}
\thanks{}

\author{Tristan Pryer} 
\address{ Tristan Pryer\newline
  School of Mathematics, Statistics \& Actuarial Science\newline 
  University of Kent\newline
  Canterbury\newline 
  GB-CT2 7NF, England UK} 
\curraddr{}
\email{\linkedemail{T.Pryer@kent.ac.uk}} 
\thanks{T.P. was also partially supported by the EPSRC grant EP/H024018/1}

\title[Energy Consistent DG methods for the NSK system]{Energy Consistent DG Methods for the Navier--Stokes--Korteweg system} 

\date{\today}
\pdfformat{true}
\begin{document}

\begin{abstract}
  We design consistent discontinuous Galerkin finite element schemes
  for the approximation of the \EK and the \NSK systems. We show that
  the scheme for the \EK system is energy and mass conservative and
  that the scheme for the \NSK system is mass conservative and
  monotonically energy dissipative. In this case the dissipation is
  isolated to viscous effects, that is, there is no numerical
  dissipation. In this sense the methods is consistent with the energy dissipation of the
  continuous PDE systems.
\end{abstract}
\maketitle


  




\section{Introduction}
\label{sec:introduction}

In this work we propose a new class of finite element methods for the
Navier-Stokes-Korteweg system which are \emph{by design consistent
  with the energy dissipation structure of the problem.} The methods
are of arbitrary high order of accuracy and provide physically
relevant approximations free of numerical artifacts.  It seems that
these are the first methods in the literature enjoying these
properties.

Liquid vapour flow occurs in many technical applications and natural
phenomena. A particularly interesting and challenging case is when the
fluid undergoes phase transition, \ie there is mass transfer between
the phases, which is driven by thermodynamics.  The applications of
these phenomena are extremely varied, for example, it is applicable to
modelling the fuel injection system in modern car engines and also to
the study of cloud formation.  The modelling of these phenomena can be
traced back to \cite{vdW88,Kor01}, however there remain open
questions, for example, what is the correct model for the given
application at hand.

The compressible flow of a single substance containing both a liquid
and vapour phase undergoing a phase transition can be modelled by
different techniques.  One widely used approach for the treatment of
these problems, which emerged in the last few decades, see
\cite{AMW98} and references therein, is the so called \emph{diffuse
  interface approach}.  In this philosophy the phases are separated by
a (thin) interfacial layer across which the fields vary smoothly.  The
benefit of this approach is that there is only one set of PDEs solved
on the whole domain whose solution already includes the position of
the interfacial layer. However, these models must include a parameter
distinguishing when we are in one phase or another.  In most diffuse
interface models this is a more or less arbitrary indicator function
based on the mass or volume fraction of one of the constituents.

In this contribution we will consider the isothermal \NSK system which
is a diffuse interface model but here the mass density serves as a
phase indicator, it originates in the work of Korteweg \cite{Kor01}
and van der Walls \cite{vdW88} and was derived in modern terminology in
\cite{DS85,TN92,JLCD01}.  This model includes surface tension effects
by a third order term in the momentum balance which corresponds to a
non-local (gradient) contribution in the energy functional. Another
feature of compressible diffuse interface models is a non--monotone
constitutive relation for the pressure. This corresponds to a
non--convex local part of the energy, see equation
\eqref{eq:double-well-to-pressure}. 
\begin{figure}[h]
  \caption{
    \label{energy_pressure_pic}
    The relation between the pressure function and the double well potential.
    }
  \begin{center}
    \subfigure[][Pressure function]{
      \includegraphics[scale=\figscale, width=0.37\figwidth]
                      {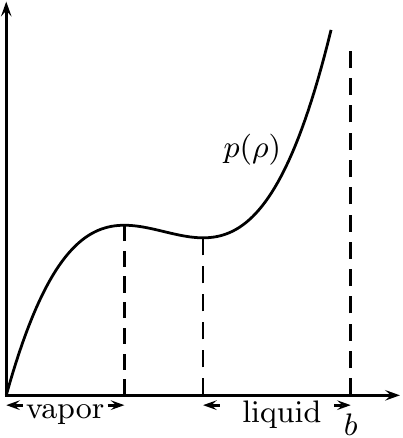}
    }
    \subfigure[][Double well potential]{
      \includegraphics[scale=\figscale, width=0.37\figwidth]
                      {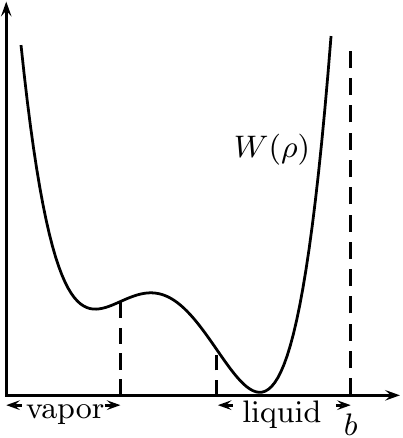}
    }
  \end{center}
\end{figure}
As can be seen from Figure \ref{energy_pressure_pic} the phases of the
problem (liquid/vapour) are the corresponding regions where the
pressure function is monotonically increasing.

The Korteweg type third order term together with the non-monotonicity
of the pressure function cause several issues in the numerical
treatment of this problem. In previous numerical studies \cite{JTB02,
  Die07, BP11} it has been observed that ``classical''
explicit-in-time finite volume (FV) and discontinuous Galerkin (DG)
schemes which  use standard fluxes used in the computational  conservation laws
introduce several numerical artifacts.

The first artifact is non-monotonicity of the energy. The \EK model is
energy conservative over time whereas the \NSK model is monotonically
energy dissipative. In the \NSK model all the dissipation is due to
viscous effects (see Lemma \ref{lem:energy}). The classical FV and DG
methods applied to the \NSK system lead to a non-monotone behaviour of
the energy. This is mainly due to the fact that these ``classical''
schemes introduce standard diffusion in the mass conservation equation
as a stabilising mechanism.  While for convex energies standard
diffusion in fact leads to energy dissipation, it may lead to an
increase in energy for multiphase flow \cite{Die07,DGR11}. In fact
standard diffusion is also present in the finite element method
proposed in \cite{BP11}.


The second artifact are so called \emph{parasitic currents}, \ie the
schemes are not well-balanced, as they do not preserve the correct
equilibria.  Parasitic currents occur when equilibrium is approached
and the numerical velocity field does not vanish uniformly, but in the
interfacial layer large velocities whose magnitude is dependent on the
gridsize and inversely dependent on the width of interfacial layer
appear \cite[\S 5]{Die07}.  As the interfacial layer is extremely thin
this effect cannot be neglected in practical computations.


Both the non-monotone behaviour of the \NSK energy and the parasitic
currents are due to numerical regularisation terms which are not
adapted to the variational structure of the problem, see \cite{DGR11}
for a study on regularisation terms taking into account the underlying
variational structure of the problem. For previous works on scalar
dispersive equations by discontinuous Galerkin methods we refer to 
 \cite{CheShu08, BCKarX2011, XuChu11}. 

The key target in the work at hand is to consider a high-order DG
discretisation of the problem which aims at preserving the energy
dissipation inequality satisfied by the original problem and avoiding
the introduction of any artificial diffusion terms. By achieving this
goal we can treat the case where the system preserves the energy
exactly. In addition, we can address the case where the system has
natural dissipation and the energy is diminishing. Our schemes are
therefore \emph{energy consistent} in the sense that they are
consistent with the energy dissipation structure of the \NSK system.
The resulting schemes are free from the above mentioned artifacts of
other approximating methods in the literature and are successful in
computing the physically relevant solution.  It is to be noted that
our approach does not hinge on an adaptation of ``entropy conservative
schemes" developed for conservation laws, \cite{Tad03}.  The
non-monotone pressure function makes a direct application of this
approach unfeasible in our case.  Conservative DG schemes for the scalar generalized KdV 
equation were suggested recently in \cite{BCKarX2011}.  To achieve our goals we follow a
constructive step-by-step approach. Motivated by the proof of energy
conservation at the continuous level we introduce a new mixed
formulation for the \NSK system. This mixed formulation will be the
basis of our discrete schemes. We first discretise in space by
employing a DG approach with generic discrete fluxes. Then we
specifically identify the properties and thus the fluxes which yield
\emph{energy consistent} schemes.
  Then we consider Crank-Nicolson type
time discretisation and identify the precise time discrete method
which is {energy consistent.}  By combining the ideas of space and
time discretization we obtain the fully discrete schemes with the
desired properties.

The structure of the paper is
as follows: In \S\ref{sec:model-prob} we introduce the \NSK model
problem as well as some of its conservative properties. We give the
mixed formulation and necessary notation which will be used throughout
the paper. In addition we describe the mixed formulation in the
\emph{broken Sobolev} framework necessary for the construction of the
DG scheme. In \S\ref{sec:numerical-methods} we detail the construction
of the energy consistent DG scheme initially in the spatially
semidiscrete case. We then move on to the temporal semidiscrete case
in \S\ref{sec:num-methods-temporal} and combine the results to obtain
an energy consistent fully discrete scheme in
\S\ref{sec:num-methods-fully}.  In \S\ref{sec:numerics} we perform
various numerical experiments to test the convergence, conservativity
and computational properties of the scheme.


\section{Model problem, mixed formulation and discretisation}
\label{sec:model-prob}

In this section we formulate the model problem, fix notation and give
some basic assumptions. Let $\W\subset \reals^{d}$, with $d=1,2,3$ be a bounded
domain. We then begin by introducing the Sobolev spaces
\cite{Ciarlet:1978,Evans:1998}
\begin{equation}
    \sobh{k}(\W) 
    := 
    \ensemble{\phi\in\leb{2}(\W)}
             {\D^{\vec\alpha}\phi\in\leb{2}(\W), \text{ for } \norm{\geovec\alpha}\leq k},
  \end{equation}
  which are equipped with norms and semi-norms
  \begin{gather}
    \Norm{u}_{k}^2
    := 
    \Norm{u}_{\sobh{k}(\W)}^2 
    = 
    \sum_{\norm{\vec \alpha}\leq k}\Norm{\D^{\vec \alpha} u}_{\leb{2}(\W)}^2 
    \\
    \AND \norm{u}_{k}^2 
    :=
    \norm{u}_{\sobh{k}(\W)}^2 
    =
    \sum_{\norm{\vec \alpha} = k}\Norm{\D^{\vec \alpha} u}_{\leb{2}(\W)}^2
\end{gather}
respectively, where $\vec\alpha = \{ \alpha_1,...,\alpha_d\}$ is a
multi-index, $\norm{\vec\alpha} = \sum_{i=1}^d\alpha_i$ and
derivatives $\D^{\vec\alpha}$ are understood in a weak sense. In addition, let
\begin{equation}
  \hoz 
  := 
  \ensemble{\phi\in\sobh1(\W)}
           {\phi\vert_{\partial\W} = 0} 
  \AND
  \hon(\W)
  := 
  \ensemble{\geovec \phi\in\qb{\sobh1(\W)}^d}
           {\Transpose{\qp{\geovec \phi \vert_{\partial\W}}}\geovec n = 0} 
\end{equation}
where $\geovec n$ denotes the outward pointing normal to $\partial
\W$.

We use the convention that for a multivariate function, $u$, the
quantity $\nabla u$ is a column vector consisting of first order
partial derivatives with respect to the spatial coordinates. The
divergence operator, $\div{}$, acts on a vector valued multivariate
function and $\Delta u := \div\qp{\nabla u}$ is the generalised
Laplacian operator. We also note that when the Laplacian acts on a
vector valued multivariate function, it is meant componentwise. 
Moreover, for a vector field $\geovec v$, we denote its Jacobian by 
$\D \geovec v$. We
also make use of the following notation for time dependant Sobolev spaces:
\begin{equation}
  \leb{2}(0,T; \sobh{k}(\W))
  :=
  \ensemble{u : [0,T] \to \sobh{k}(\W)}
           {\int_0^T \Norm{u(t)}_{k}^2 \d t < \infty}.
\end{equation}

\subsection{Model problem}

Consider a fluid in the domain $\W$ with density $\rho$ and
velocity $\geovec v$.  The \NSK system is made up of the balances of
mass and momentum of said fluid, that is,
\begin{equation}
  \label{eq:NSK}
  \begin{split}
  \pdt{\rho} 
  +
  \div\qp{\rho \geovec v} 
  &= 0
  \\
  \pdt{\qp{\rho \geovec v}} 
  +
  \div\qp{\rho \geovec v\otimes \geovec v}
  +
  \nabla p(\rho) 
  &=
  \mu \Delta \geovec v 
  +
  \gamma \rho \nabla \Delta \rho
  \end{split} 
  \quad 
  \text{ in }
  \quad \W \times (0,T)
\end{equation}
where $p$ is a non--monotone pressure function, $\mu$ is a viscosity
coefficient and $\gamma$ a capillarity coefficient.  The pressure
function $p$ is linked to a double well-potential $W=W(\rho)$ via the relation (Figure \ref{energy_pressure_pic})
\begin{equation}
  \label{eq:double-well-to-pressure}
  p(\rho)
  =
  \rho W'(\rho) - W(\rho).
\end{equation}

Let $\geovec n$ be the outward pointing normal to $\partial\W$, suppose
the system \eqref{eq:NSK} is given with boundary conditions
\begin{equation}
  \label{eq:boundary conditions_1}
  \geovec v 
  =
  \geovec 0
  \AND
  \Transpose{\qp{\nabla \rho}}\geovec n
  =
  0
  \text{ on }
  \partial\W \times (0,T)
\end{equation}
and initial conditions
\begin{equation}
  \label{eq:initial conditions_1}
  \rho(\cdot,0)=\rho^0,\quad \geovec v(\cdot,0)=\geovec v^0
  \text{ in }
  \W 
\end{equation}
for given functions $\rho^0 \in \sobh{1}(\W)$ and $\geovec v^0 \in
(\sobh{1}(\W))^d$ such that $W(\rho^0)\in\leb{1}(\W)$. The system
(\ref{eq:NSK}) conserves mass as well as satisfying a momentum balance
together with an energy dissipation equality, \ie
\begin{gather}
  \ddt 
  \qp{\int_\W \rho \d \geovec x}
  =
  0 
  \\
  \ddt
  \qp{\int_\W \rho \geovec v \d \geovec x}
  =
  \mu \int_{\partial \W} (\D\geovec v) \geovec n \d s 
  \\
\label{eq:energy_dissipation_equality}
  \ddt
  \qp{\int_\W W(\rho) 
  + 
  \frac{1}{2} 
  \rho \norm{\geovec v}^2 
  + 
  \frac{\gamma}{2} 
 \norm{\nabla \rho}^2 \d \geovec x}
  =
  - \mu \int_{\W} \norm{\D \geovec v}^2 \d \geovec x
  ,
\end{gather}
respectively.  The energy dissipation equality is only valid for
smooth solutions. In case the system permits shocks, they would
trigger additional energy dissipation and
\eqref{eq:energy_dissipation_equality} would have to be replaced by an
inequality. While the first two equalities follow by integrating the
mass and momentum balance \eqref{eq:NSK}. The derivation of the energy
dissipation equality is a little bit more involved. For completeness
the result is formulated as Lemma \ref{lem:energy}.  Moreover, the
proof of Lemma \ref{lem:energy} serves as a guideline for the
construction of energy consistent discrete schemes.

\begin{Hyp}[finite Helmholtz energy]
  From hereon in we will assume that for a given $\rho$ we have that
  $W(\rho) \in \leb{1}(0,T;\leb{1}(\W))$.
\end{Hyp}

\begin{Lem}
  \label{lem:energy}
  For every smooth solution $(\rho,\geovec v) \in
  \leb2(0,T;\sobh3(\W)) \times \leb2(0,T;\sobh2(\W))^d$ of
  \eqref{eq:NSK} such that $(\pdt \rho,\partial_t \geovec v) \in
  \leb2(0,T;\leb2(\W)) \times \leb2(0,T;\leb2(\W))^d$ which satisfies
  the boundary conditions \eqref{eq:boundary conditions_1} we have
  \begin{equation}
    \ddt 
    \qp{
      \int_\W W(\rho) 
      +
      \frac{1}{2}\rho \norm{\geovec v}^2 
      +
      \frac{\gamma}{2} \norm{\nabla \rho}^2 \d \geovec x
    }
    =
    - \mu \int_{\W} 
    \norm{\D \geovec v}^2 \d \geovec x.
  \end{equation}
\end{Lem}

\begin{Proof}
  Let us first note that the second equation of \eqref{eq:NSK} can be
  reformulated as
  \begin{equation}
    \label{eq:energy_proof_1}
    \rho \partial_t \geovec v
    + 
    \div(\rho \geovec v \otimes \geovec v) 
    -
    \div(\rho \geovec v) \geovec v
    +
    \rho \nabla W'(\rho) 
    -
    \mu \Delta \geovec v 
    -
    \gamma \rho \nabla \Delta \rho
    =
    0.
  \end{equation}
  Multiplying the first equation of \eqref{eq:NSK} by 
  $W'(\rho) + \frac{1}{2}  \norm{\geovec v}^2 - \gamma \Delta \rho$
  we see
  \begin{equation}
    \label{eq:energy_proof_4}
    \begin{split}
    0
    &=  
    W'(\rho) \pdt{\rho}
    +
    W'(\rho) \div\qp{\rho \geovec v}
    +
    \frac{1}{2} \norm{\geovec v}^2 \pdt{\rho} 
    +
    \frac{1}{2}  \norm{\geovec v}^2 \div\qp{\rho \geovec v}
    \\
    & \qquad -
    \gamma \Delta \rho \pdt{\rho}
    -
    \gamma \Delta \rho \div\qp{\rho \geovec v}
    .
    \end{split}
  \end{equation}
  Then by multiplying \eqref{eq:energy_proof_1} by $\geovec v$ and
  summing together with (\ref{eq:energy_proof_4}) we obtain
  \begin{equation}
    \label{eq:energy_proof_2}
    \begin{split}
      0
      &=
      W'(\rho)\pdt \rho 
      +
      \frac{1}{2}\pdt \rho  \norm{\geovec v}^2 
      -
      \gamma\pdt \rho  \Delta \rho
      +
      W'(\rho)\div(\rho \geovec v)
      -
      \frac{1}{2} \div(\rho \geovec v)  \norm{\geovec v}^2 
      \\
      &\qquad -
      \gamma \div(\rho \geovec v) \Delta \rho
      +
      \rho \Transpose{\geovec v} \partial_t \geovec v 
      +
      \Transpose{ \geovec v }\div(\rho \geovec v \otimes \geovec v)
      + 
      \rho \Transpose{\geovec v} \nabla W'(\rho)
      \\
      &\qquad-
      \mu \Transpose{\geovec v} \Delta \geovec v
      -
      \gamma \rho \Transpose{\geovec v} \nabla \Delta \rho.
    \end{split}
  \end{equation}
  We integrate \eqref{eq:energy_proof_2} over $\W$
  and by Greens formula we have
  \begin{equation}
    \label{eq:energy_proof_3}
    \begin{split}
      0
      &=
      \int_\W 
      W'(\rho)\pdt \rho
      +
      \frac{1}{2}\pdt \rho \norm{\geovec v}^2 
      +
      \rho \Transpose{\geovec v} \partial_t \geovec v 
      +
      \frac{\gamma}{2}\partial_t \norm{\nabla \rho}^2
      +
      \mu \norm{ \D \geovec v}^2  \d \geovec x
      \\
      &    
      + 
      \int_{\partial \W}
    \qp{ \Transpose{ \qp{
      - \gamma  \pdt \rho \nabla \rho  
        +
        \frac{1}{2} \rho \geovec v \norm{\geovec v}^2
        -
        \gamma \rho \geovec v \Delta \rho 
        +
        \rho \geovec v W'(\rho)}}
        - 
      \mu  \Transpose{\geovec v} \D{ \geovec v }
      }
      \geovec n \d s.
 \end{split}
\end{equation}
The boundary integral in \eqref{eq:energy_proof_3} vanishes because of the boundary conditions.
\end{Proof}

\begin{Rem}[stable steady states]
  The energy dissipation equality gives rise to the fact that the
  (stable) steady states of \eqref{eq:NSK} are minimizers of the
  energy functional
  \begin{equation}
    \label{eq:def_energy} 
    E[\rho,\geovec v] 
    :=
    \int_\W W(\rho) 
    +
    \frac{1}{2}\rho \norm{\geovec v}^2 
    +
    \frac{\gamma}{2} \norm{\nabla \rho}^2 \d \geovec x,
  \end{equation}
under the constraint
  \begin{equation}
    \label{eq:def_energy_constr} 
    \int_\W \rho  \d \geovec x =m,
  \end{equation}
for some given $m>0$
  and therefore satisfy the Euler-Lagrange equations
  \begin{gather}
    \geovec v
    =
    \geovec 0
    \\
    \label{eq:lag-mult-EL}
    {W'(\rho) - \gamma \Delta \rho} 
    =
    \lambda,
  \end{gather}
  where $\lambda$ is the Lagrange multiplier associated with the mass
  conservation contraint (\ref{eq:def_energy_constr}).

  Note that (\ref{eq:lag-mult-EL}) is equivalent to
  \begin{equation}
    \begin{split}
      \geovec 0 
      &= 
      \nabla \qp{    {W'(\rho) - \gamma \Delta \rho} }
      \\
      &=
      \nabla p(\rho) - \gamma \rho\nabla\Delta\rho
    \end{split}
  \end{equation}
  using the relation
  \begin{equation}
    \label{eq:pressure-W-rel}
    \nabla p(\rho) = \rho \nabla W'(\rho)
  \end{equation}
  which is readily derived from (\ref{eq:double-well-to-pressure}).
\end{Rem}
  
\subsection{Classical solvability of the problem}

The well--posedness of the \NSK system and similar systems was
considered by several authors \cite{BDDJ07, BDL03, DD01, Fei02, HL96,
  Kot06}. For completeness we will state some results.

\begin{The}[existence of a solution to the \EK system \cite{BDDJ07}]
  \label{the:exist-EK}
Let $s>\tfrac{d}{2}+1$ and
  \begin{equation}
    H_s 
    := 
    H^{s+1} (\rR^d) \times H^s (\rR^d,\rR^d).
  \end{equation}
  Suppose the initial
  data $(\rho_0,\geovec v_0) \in (\underline \rho(0),\geovec
  {\underline v}(0)) + H_{s}$ where $\underline \rho,
  \geovec{\underline v}$ is a special solution such that $\underline
  \rho$ is bounded away from zero and the Hessian of $\underline\rho$,
  $ \D\nabla \underline \rho$, as well as the Jacobian of
  $\geovec{\underline v}$, $\D \geovec{\underline v}$ are both 
  $C([0,T],H^{s+3}(\rR^d,\rR^{d\times d}))$ for some $T>0.$ Then the
  \EK system admits a unique solution $(\rho,\geovec v) \in
  (\underline \rho,\geovec {\underline v})+ C^1([0,T),H_{s-2}) \cap
    C([0,T),H_s)$ satisfying the initial data $(\rho_0,\geovec v_0)$.
\end{The}

\begin{The}[existence of a solution to the \NSK system \cite{DD01}]
  \label{the:exist-NSK}
  Let $B^s =B^s_{2,1}(\rR^d)$ denote the homogeneous Besov space. Let
  $\bar \rho >0$ be a reference density such that $p'(\bar \rho)>0$.
  Suppose also
  that the initial data $\rho_0,\geovec v_0$ satisfies $\rho_0 - \bar
  \rho \in B^{d/2},$ $\rho_0 \geq c >0$ and $\geovec v_0 \in
  (B^{d/2-1})^d$.

  Then there exists a $T>0$ such that the \NSK system has a unique
  solution $(\rho,\geovec v)$ with initial data $\rho_0,\geovec v_0$
  such that $\rho - \bar \rho \in C([0,T),B^{d/2}) \cap L^
    1([0,T),B^{d/2+2})$ and $\geovec v \in C([0,T),(B^{d/2-1}))^d \cap
        L^1([0,T),(B^{d/2+1}))^d$.
\end{The}

\begin{Rem}[]
  Theorems \ref{the:exist-EK} and \ref{the:exist-NSK} motivate us to
  construct numerical schemes which are adapted to the smooth
  situation. In particular, enforcing the energy dissipation equality
  proven in Lemma \ref{lem:energy}.
\end{Rem}

\subsection{Mixed formulation}
To mimic the proof of Lemma \ref{lem:energy} at the discrete level
 it will be essential to have at our disposal  a numerical formulation in which $\geovec v$ and $\tau$,
  which depend nonlinearly on the original variables $\rho$ and $\rho
  \geovec v$, are permitted as test functions.  Indeed, this is our main motivation to reformulate
the \NSK system \eqref{eq:NSK}    as a mixed system
of PDEs by the introduction of two auxilliary variables, $\tau \AND
\geovec q$, and using the relation of the pressure function and the
double well potential (\ref{eq:pressure-W-rel}).

The mixed formulation is then to seek $\qp{\rho, \geovec v, \tau, \geovec q}$ such that
\begin{gather}
  \label{eq:mixed-system-1}
  \pdt{\rho} 
  +
  \div\qp{\rho \geovec v}
  =
  0
  \\
  \label{eq:mixed-system-2}
  \rho 
  \pdt{\geovec v}
  +
  \div\qp{\rho \geovec v\otimes \geovec v}
  -
  \div\qp{\rho \geovec v} \geovec v
  +
  \rho\nabla \tau
  -
  \frac{1}{2}
  \rho \nabla \norm{\geovec v}^2 - \mu \Delta \geovec v 
  =
  0
  \\
  \label{eq:mixed-system-3}
  \tau 
  -
  W'(\rho)
  +
  \gamma \div\qp{\geovec q}
  -
  \frac{1}{2}\norm{\geovec v}^2
  =
  0
  \\
  \label{eq:mixed-system-4}
  \geovec q
  -
  \nabla \rho
  =
  0
\end{gather}
which is coupled with the boundary conditions
\begin{equation}
  \geovec v 
  =
  \geovec 0
  \AND
  \Transpose{\geovec q} \geovec n 
  =
  0
  \text{ on }
  \partial\W \times (0,T)
\end{equation}
and the initial conditions \eqref{eq:initial conditions_1}.

%

\begin{Rem}[alternate notation]
  We note that the second and third term on the left hand side of
  (\ref{eq:mixed-system-2}) can be rewritten as
  \begin{equation}
    \div\qp{\rho \geovec v\otimes \geovec v}
    -
    \div\qp{\rho \geovec v} \geovec v
    =
    \rho \qp{\Transpose{\geovec v}\nabla} \geovec v,
  \end{equation}
  which is the standard notation in the incompressible scenario.
\end{Rem}

\subsection{Discretisation}

Let $\T{}$ be a conforming, shape regular triangulation of $\W$,
namely, $\T{}$ is a finite family of sets such that
\begin{enumerate}
\item $K\in\T{}$ implies $K$ is an open simplex (segment for $d=1$,
  triangle for $d=2$, tetrahedron for $d=3$),
\item for any $K,J\in\T{}$ we have that $\closure K\meet\closure J$ is
  a full subsimplex (i.e., it is either $\emptyset$, a vertex, an
  edge, a face, or the whole of $\closure K$ and $\closure J$) of both
  $\closure K$ and $\closure J$ and
\item $\union{K\in\T{}}\closure K=\closure\W$.
\end{enumerate}
We use the convention where $\funk h\W\reals$ denotes the
\emph{meshsize function} of $\T{}$, \ie 
\begin{equation}
  h(\geovec{x}):=\max_{\closure K\ni \geovec x}h_K,  
\end{equation}
where $h_K$ is the diameter of an element $K$. We let $\E{}$ be the
skeleton (set of common interfaces) of the triangulation $\T{}$ and
say $e\in\E$ if $e$ is on the interior of $\W$ and $e\in\partial\W$ if
$e$ lies on the boundary $\partial\W$.

\begin{Defn}[Broken Sobolev spaces, trace spaces]
  \label{defn:broken-sobolev-space}
  We introduce the
  broken Sobolev space
  \begin{equation}
    \sobh{k}(\T{})
    :=
    \ensemble{\phi}
             {\phi|_K\in\sobh{k}(K), \text{ for each } K \in \T{}},
  \end{equation}
  similarly for $\sobh1_0(\T{})$ and $\hon(\T{})$.

  We also make use of functions defined in these broken spaces
  restricted to the skeleton of the triagulation. This requires an
  appropriate trace space
  \begin{equation}
    \Tr{\E}
    :=
    \prod_{K\in\T{}} \leb{2}(\partial K) \subset \prod_{K\in\T{}} \sobh{\frac{1}{2}}(\partial K).
  \end{equation}
\end{Defn}

Let $\poly p(\T{})$ denote the space of piecewise polynomials of
degree $p$ over the triangulation $\T{}$ we then introduce the
\emph{finite element spaces}
\begin{gather}
  \label{eqn:def:finite-element-space}
  \fes := \dg{p} = \poly p(\T{})
  \\
  \feszero := \fes \cap \hoz(\T{})
  \\
  \fesn := \fes^d \cap \hon(\T{})
\end{gather}
to be the usual spaces of (discontinuous) piecewise polynomial
functions. For simplicity we will assume that $\fes$ is constant in time.

\begin{Defn}[jumps and averages]
  \label{defn:averages-and-jumps}
  We may define average and jump operators over $\Tr{\E}$ for
  arbitrary scalar, $v\in\Tr{\E}$, and vector valued functions, $\geovec
  v\in\Tr{\E}^d$.
  \begin{equation}
    \label{eqn:average}
    \dfunkmapsto[.]
	        {\avg{\cdot}}
	        v
	        {\Tr{\E}}
	        {\frac{1}{2}\qp{v|_{K_1} + v|_{K_2}}}
	        {\leb{2}(\E)}
  \end{equation}
  
  \begin{equation}
    \label{eqn:average-vec}
    \dfunkmapsto[.]
	        {\avg{\cdot}}
	        {\geovec v}
	        {\qp{\Tr{\E}}^d}
	        {\frac{1}{2}\qp{\geovec{v}|_{K_1} + \geovec{v}|_{K_2}}}
	        {\qp{\leb{2}(\E)}^d}
  \end{equation}
  
  \begin{equation}
    \label{eqn:jump}
    \dfunkmapsto[.]
	        {\jump{\cdot}}
	        {v}
	        {\qp{\Tr{\E}}}
	        {{{v}|_{K_1} \geovec n_{K_1} + {v}|_{K_2}} \geovec n_{K_2}}
	        {\qp{\leb{2}(\E)}^d}
  \end{equation}
  
  \begin{equation}
    \label{eqn:jump-vec}
    \dfunkmapsto[.]
	        {\jump{\cdot}}
	        {\geovec v}
	        {\qp{\Tr{\E}}^d}
	        {{\Transpose{\qp{\geovec{v}|_{K_1}}}\geovec n_{K_1} + \Transpose{\qp{\geovec{v}|_{K_2}}}\geovec n_{K_2}}}
	        {\qp{\leb{2}(\E)}}
  \end{equation}
  
  \begin{equation}
    \label{eqn:tensor-jump}
    \dfunkmapsto[,]
	        {\tjump{\cdot}}
	        {\vec v}
	        {\qp{\Tr{\E\cup \partial\W}}^d}
	        {{\vec{v}|_{K_1} }\otimes \geovec n_{K_1} + \vec{v}|_{K_2} \otimes\geovec n_{K_2}}
	        {\qp{\leb{2}(\E)}^{d\times d}}
  \end{equation}
  where $\geovec n_{K_i}$ denotes the outward pointing normal to $K_i$.
  Note that on the boundary of the domain $\partial\W$ the jump and
  average operators are defined as
  \begin{gather}
    \jump{v}
    \Big\vert_{\partial\W}
    := 
    v\geovec n
    \qquad 
    \jump{\geovec v}
    \Big\vert_{\partial\W} 
    :=
    \Transpose{\geovec v}\geovec n
    \\
    \avg{v}
    \Big\vert_{\partial\W} 
    := v
    \qquad 
    \avg{\geovec v}
    \Big\vert_{\partial\W}
    :=
    \geovec v.
  \end{gather}
\end{Defn}

\subsection{Elementwise formulation and discrete fluxes}

As a next step towards the construction of a numerical scheme we give the elementwise variational formulation to the problem in mixed form
(\ref{eq:mixed-system-1})--~(\ref{eq:mixed-system-4}). It requires to find
$\qp{\rho, \geovec v} \in \leb{2}(0,T;\sobh{1}(\T{})) \times
\qp{\leb{2}(0,T;\hoz(\T{})}^d$ with $\qp{\pdt{\rho}, \pdt{\geovec v}}
\in \leb{2}(0,T;\leb{2}(\T{})) \times \qp{\leb{2}(0,T;\leb2(\T{})}^d$
and $\qp{\tau, \geovec q} \in \leb{2}(0,T;\sobh{1}(\T{}))\times {\leb{2}(0,T;\hon(\T{}))}$ such
that $W'(\rho) \in \leb{2}(0,T;\leb{2}(\T{}))$ and

\begin{equation}
  \begin{split}
    \label{eq:variational_formulation}
    0 
    &=
    \int_\W \qp{{\pdt\rho} + {\div\qp{\rho \geovec v}}} {\psi} \d \geovec x
    + 
    \int_\E
    F_1\qp{\rho, \geovec v, \tau, \geovec q, \psi} \d s
    \Foreach \psi \in \sobh{1}(\T{})
    \\
    0
    &=
    \int_\W
    \Transpose{\qp{\rho \pdt{{\geovec v}}
      +
      \div\qp{\rho \geovec v\otimes \geovec v}
      -
      \div\qp{\rho \geovec v} \geovec v
      +
      \rho\nabla \tau
      -
      \frac{1}{2}
      \rho \nabla \norm{\geovec v}^2 
    }} \geovec \chi
    \d \geovec x
    \\
    & \qquad \qquad + 
    \int_\E 
    F_2\qp{\rho, \geovec v, \tau, \geovec q, \geovec \chi} \d s  + \mu B(\geovec v,\geovec \chi)
    \Foreach \geovec \chi \in \qp{\hoz(\T{})}^d
    \\
    0 
    &=
    \int_\W
    \qp{  \tau - W'(\rho) + \gamma \div\qp{\geovec q} - \frac{1}{2}\norm{\geovec v}^2 } \xi \d \geovec x
    \\
    & \qquad\qquad+
    \int_\E
    F_3\qp{\rho, \geovec v, \tau, \geovec q, \xi} \d s 
    \Foreach \xi \in \sobh{1}(\T{})
    \\
    0 
    &=
    \int_\W \Transpose{\qp{\geovec q - \nabla \rho}} \geovec \zeta \d \geovec x
    +
    \int_\E  F_4\qp{\rho, \geovec v, \tau, \geovec q, \geovec \zeta} \d s
    \Foreach \geovec \zeta \in {\hon(\T{})}
  \end{split}
\end{equation}
where 
\begin{equation}
  \begin{split}
    F_1, F_3 &: \sobh{1}(\T{})\times \hoz(\T{})^d \times \sobh{1}(\T{}) \times \hon(\T{}) \times \sobh{1}(\T{}) \to \leb2(\E)\\
    F_2 &: \sobh{1}(\T{})\times \hoz(\T{})^d \times \sobh{1}(\T{}) \times \hon(\T{}) \times \hoz(\T{})^d  \to \leb2(\E)
    \\
    F_4 &: \sobh{1}(\T{})\times \hoz(\T{})^d \times \sobh{1}(\T{}) \times \hon(\T{}) \times \hon(\T{})  \to \leb2(\E)
  \end{split}
\end{equation}
are appropriate choices of elementwise fluxes to be chosen in the
sequel to suit our purposes; the operators $\div$ and $\nabla$ are
understood henceforth to be defined elementwise and $B:
\qp{\hoz(\T{})}^d \times \qp{\hoz(\T{})}^d \to \rR$ is a bilinear
form, corresponding to a weak formulation of the Laplacian.  We will
also assume that the fluxes $F_1,\dots,F_4$ only depend on the traces
of their arguments and are linear in the test functions.
We would like to mention that  the spaces for the variational formulation are chosen such that all integrals in 
\eqref{eq:variational_formulation} are well-defined. We do not claim that there is a well-posedness analysis for \eqref{eq:variational_formulation}
with the given spaces; \eqref{eq:variational_formulation} serves only as a basis to define the spatial discrete DG scheme in the next section.


\section{Development of  energy consistent numerical methods -- the spatially discrete case}
\label{sec:numerical-methods}

In this section we will detail the methodology behind the construction
of the energy consistent finite element scheme. We present our main
results which show the conditions a generic scheme applied to the
variational formulation (\ref{eq:variational_formulation}) with no
diffusion (i.e., $\mu = 0$) is mass and energy conservative.  If a
scheme conserves mass for $\mu=0$ this does not change for $\mu
\not=0$, as the mass conservation equation is not affected by a
reasonable discretization of the viscosity. Moreover if a scheme
conserves energy for $\mu=0$, for $\mu \not=0$ all energy dissipation
is due to viscosity.  For simplicity we first detail the calculations
for the spatially discrete case, then construct a temporally discrete
scheme. We also give a condition when a scheme falling under our
framework can also conserve momentum.

\subsection{Spatially discrete scheme}

Throughout the calculations in this section we will regularly refer to
the following proposition.
\begin{Pro}[elementwise integration]
  \label{Pro:trace-jump-avg}
  Let 
  \begin{equation}
  \sobh{\div}(\T{}) 
  := 
  \ensemble{\geovec p \in(\leb{2}(\T{}))^d}{\div{\geovec p} \in \leb{2}(\T{})}.
  \end{equation}
  Suppose $\geovec p 
  \in \sobh{\div}(\T{})$
  and $\phi\in\sobh1(\T{})$ then 
  \begin{equation}
    \begin{split}
      \sum_{K\in\T{}}
      \int_K \div\qp{\geovec p} \phi \d \geovec x
      =
      \sum_{K\in\T{}}
      \qp{
        -
        \int_K
        \Transpose{\geovec p} \nabla \phi \d \geovec x
        +
        \int_{\partial K}
        \phi \Transpose{\geovec p} \geovec n_K \d s
        }.
    \end{split}
  \end{equation}
  In particular we have $\geovec p \in \Tr{\E}^d$ and $\phi \in
  \Tr{\E}$, and the following identity holds
  \begin{equation}
    \sum_{K\in\T{}}
    \int_{\partial K}
    \phi \Transpose{{\geovec p}}
    \geovec n_K \d s
    =
    \int_\E 
    \jump{\geovec p} 
    \avg{\phi}
    \d s
    +
    \int_{\E\cup\partial\W}
    \Transpose{\jump{\phi}}
    \avg{\geovec p}
    \d s
    =
    \int_{\E\cup\partial\W}
    \jump{\geovec p \phi}
    \d s.
  \end{equation}
\end{Pro}

\subsection{General numerical scheme}
A generic spatially discrete DG formulation to the problem in mixed
form (\ref{eq:mixed-system-1})--(\ref{eq:mixed-system-4}) is to find
$\rho_h, \tau_h:[0,T] \to \fes$ and $\geovec v_h : [0,T]
\to \feszero^d \AND \geovec q_h : [0,T] \to \fesn$ such that
\begin{equation}
  \begin{split}
    \label{eq:discrete_family}
    0 
    &=
    \int_\W \qp{{\pdt\rho_h} + {\div\qp{\rho_h \geovec v_h}}} {\Psi} \d \geovec x
    + 
    \int_\E
    F_1\qp{\rho_h, \geovec v_h, \tau_h, \geovec q_h, \Psi} \d s
    \Foreach \Psi \in \fes
    \\
    0
    &=
    \int_\W
    \Transpose{\qp{\rho_h \pdt{{\geovec v_h}}
      +
      \div\qp{\rho_h \geovec v_h\otimes \geovec v_h}
      -
      \div\qp{\rho_h \geovec v_h} \geovec v_h
      +
      \rho_h\nabla \tau_h
      -
      \frac{1}{2}
      \rho_h \nabla \norm{\geovec v_h}^2 
    }} \geovec \Chi
    \d \geovec x
    \\
    & \qquad \qquad + 
    \int_\E 
    F_2\qp{\rho_h, \geovec v_h, \tau_h, \geovec q_h, \geovec \Chi} \d s
    + \mu B_h (\geovec v_h,\geovec \Chi)
    \Foreach \geovec \Chi \in \feszero^d
    \\
    0 
    &=
    \int_\W
    \qp{  \tau_h - W'(\rho_h) + \gamma \div\qp{\geovec q_h} - \frac{1}{2}\norm{\geovec v_h}^2 } \Xi \d \geovec x
    \\
    &\qquad\qquad +
    \int_\E
    F_3\qp{\rho_h, \geovec v_h, \tau_h, \geovec q_h, \Xi} \d s 
    \Foreach \Xi \in \fes
    \\
    0 
    &=
    \int_\W \Transpose{\qp{\geovec q_h - \nabla \rho_h}} \geovec \Zeta \d \geovec x
    +
    \int_\E  F_4\qp{\rho_h, \geovec v_h, \tau_h, \geovec q_h, \geovec \Zeta} \d s
    \Foreach \geovec \Zeta \in \fesn,
  \end{split}
\end{equation}
where $B_h: \feszero^d \times \feszero^d \to \rR$ is a discretization of $B$.

\subsection{Consistency and conservation}

Now we will give abstract properties of the fluxes which determine
whether the scheme is consistent and conserves mass, momentum or
energy.

\begin{Defn}[consistency]
  \label{The:consistency}
  A generic scheme having the form (\ref{eq:discrete_family}) is said to be 
  consistent provided
  \begin{equation}
    \label{eq:cond_consistency}
    F_i(\rho,\geovec v, \tau, \geovec q, \cdot  ) 
    \equiv
    0 \text{ for } i = 1,\dots,4
  \end{equation}
 for all smooth functions $\rho,\tau \in\leb{2}(0,T; \sobh1(\W))$ and
 $\geovec v,\geovec q \in \qb{\leb{2}(0,T; \sobh1(\W))}^d$ .
\end{Defn}

\begin{The}[conservation]
  \label{The:conservation}
 For $\mu=0$ a generic scheme of the form (\ref{eq:discrete_family}) conserves:
  \begin{enumerate}
  \item 
    Mass, that is,
  \begin{equation}
    \label{eq:mass-cons}
    \ddt \qp{\int_\W \rho_h\d \geovec x} = 0 
  \end{equation}
  if and only if
  \begin{equation}
    \label{eq:cond_mass}
    \int_\E 
    F_1(\rho_h,\geovec v_h, \tau_h, \geovec q_h, 1) 
    \d s
    =
    -
    \int_\E
    \jump{\rho_h \geovec v_h}
    \d s
    \Foreach \rho_h,\tau_h \in \fes, \geovec v_h \in \feszero^d, \geovec q_h \in \fesn,
  \end{equation}
  where $1$ is the the constant element of $\sobh1(\T{})$ which is $1$
  everywhere.
\item Energy, that is,
  \begin{equation}
    \label{eq:energy-cons}
    \ddt \qp{\int_\W W(\rho_h) 
      + \frac{1}{2}\rho_h \norm{\geovec v_h}^2 
      + \frac{\gamma}{2} \norm{\geovec q_h}^2 \d \geovec x }= 0
  \end{equation}
  if and only if
  \begin{equation}
    \label{eq:energy_cons}
    \begin{split}
      \int_\E 
      F_1(\rho_h,\geovec v_h, \tau_h, \geovec q_h, \tau_h) 
      +
      F_2(\rho_h,\geovec v_h, \tau_h, \geovec q_h, \geovec v_h) 
      +
      \jump{\rho_h \tau_h \geovec v_h}\d s  
      &=
      0 
      , 
      \\
      \int_\E
      F_3(\rho_h,\geovec v_h, \tau_h, \geovec q_h, \partial_t \rho_h) 
      -
      \gamma  D_t F_4(\rho_h,\geovec v_h, \tau_h, \geovec q_h ,\geovec q_h )
      + 
      \gamma \jump{\pdt{ \rho_h} \geovec q_h} \d s
      &=
      0,
    \end{split}
  \end{equation}
  for all $\rho_h,\tau_h: [0,T] \to \fes, \geovec v_h:
  [0,T] \to \feszero^d \AND \geovec q_h : [0,T] \to \fesn$. Note that we use the notation $D_t
  F_4(\rho_h,\geovec v_h, \tau_h, \geovec q_h ,\geovec Z )$ for the
  time derivative since $\rho_h,\geovec v_h, \tau_h, \geovec q_h$ are
  time dependent but $\geovec Z$ is independent of time, as $\fesn$ is independent of time. 
\end{enumerate}
\end{The}

\begin{Cor}[Energy dissipation]\label{cor:energy-dissipation}
 Let $B_h$ be a coercive discretisation of $B$, then for $\mu>0$ a
 generic scheme of the form (\ref{eq:discrete_family}) conserves mass,
 if and only if \eqref{eq:cond_mass} is satisfied and it satisfies the
 energy dissipation equality
  \begin{equation}
    \label{eq:energy-diss}
    \ddt \qp{\int_\W W(\rho_h) 
      + \frac{1}{2}\rho_h \norm{\geovec v_h}^2 
      + \frac{\gamma}{2} \norm{\geovec q_h}^2 \d \geovec x }= -\mu B_h (\geovec v_h,\geovec v_h) \leq 0
  \end{equation}
  if and only if
  \begin{equation}
    \begin{split}
      \int_\E 
      F_1(\rho_h,\geovec v_h, \tau_h, \geovec q_h, \tau_h) 
      +
      F_2(\rho_h,\geovec v_h, \tau_h, \geovec q_h, \geovec v_h) 
      +
      \jump{\rho_h \tau_h \geovec v_h}\d s  
      &=
      0 
      , 
      \\
      \int_\E
      F_3(\rho_h,\geovec v_h, \tau_h, \geovec q_h, \partial_t \rho_h) 
      -
      \gamma  D_t F_4(\rho_h,\geovec v_h, \tau_h, \geovec q_h ,\geovec q_h )
      + 
      \gamma \jump{\pdt{ \rho_h} \geovec q_h} \d s
      &=
      0,
    \end{split}
  \end{equation}
holds.
\end{Cor}

\begin{Cor}[Consistency, conservation and dissipation]
  \label{cor:consistency-and-conservation}
  For $\mu = 0$ the following spatially discrete scheme 
  \begin{equation}
    \label{eq:discrete_variational_formulation}
    \begin{split}
      0 
      &=
      \int_\W \qp{{\pdt\rho_h} + {\div\qp{\rho_h \geovec v_h}}} {\Psi} \d \geovec x
      - 
      \int_\E
      \jump{\rho_h \geovec v_h}\avg{\Psi} \d s
      \Foreach \Psi \in \fes
      \\
      0
      &=
      \int_\W
      \Transpose{\qp{\rho_h \pdt{{\geovec v_h}}
          +
          \div\qp{\rho_h \geovec v_h\otimes \geovec v_h}
          -
          \div\qp{\rho_h \geovec v_h} \geovec v_h
          +
          \rho_h\nabla \tau_h
          -
          \frac{1}{2}
          \rho_h \nabla \norm{\geovec v_h}^2 
      }} \geovec \Chi
      \d \geovec x
      \\
      & \qquad \qquad
      -
      \int_\E 
      \Transpose{\jump{\tau_h}}\avg{\rho_h\geovec\Chi}\d s
      + \mu B_h(\geovec v_h, \geovec \Chi)
      \Foreach \geovec \Chi \in \feszero^d
      \\
      0 
      &=
      \int_\W
      \qp{  \tau_h - W'(\rho_h) + \gamma \div\qp{\geovec q_h} - \frac{1}{2}\norm{\geovec v_h}^2 } \Xi \d \geovec x
      \\
      &\qquad\qquad -
      \int_\E
      \gamma\jump{\geovec q_h}\avg{\Xi}\d s 
      \Foreach \Xi \in \fes
      \\
      0 
      &=
      \int_\W \Transpose{\qp{\geovec q_h - \nabla \rho_h}} \geovec \Zeta \d \geovec x
      +
      \int_\E 
      \Transpose{{\jump{\rho_h}}}\avg{\geovec \Zeta} \d s
      \Foreach \geovec \Zeta \in \fesn
    \end{split}
  \end{equation}
  is consistent and conserves mass (\ref{eq:mass-cons}) and energy
  (\ref{eq:energy-cons}).

  Let $B_h(u,w)$ be the symmetric interior penalty method for the
  componentwise Laplacian given by
  \begin{equation}
    \label{eq:symm-ip-method}
    B_h(\geovec u,\geovec w) 
    =
    \int_\W 
    \frob{\D \geovec u}{\D \geovec w}
    \d \geovec x
    -
    \int_{\E\cup\partial\W} 
    \frob{\avg{\D \geovec w}}{\tjump{\geovec u}}
    +
    \frob{\avg{\D \geovec u}}{\tjump{\geovec w}}
    -
    \frac{\sigma}{h}\frob{\tjump{\geovec u}}{\tjump{\geovec w}}
    \d s,
  \end{equation}
  where $\frob{}{}$ denotes the Frobenius inner product between two
  $d\times d$ matrices, \ie $\frob{\geomat X}{\geomat Y} :=
  \trace\qp{\Transpose{\geomat X}{\geomat Y}}$. It is well know for large
  (enough) $\sigma$ this is a coercive discretisation of the componentwise
  Laplacian. Thus, for $\mu > 0$, the numerical scheme
  (\ref{eq:discrete_variational_formulation}) with
  (\ref{eq:symm-ip-method}) is consistent, conserves mass
  (\ref{eq:mass-cons}) and dissipates energy (\ref{eq:energy-diss}).
\end{Cor}

\begin{Proof}[of Theorem \ref{The:conservation}]
Let us first consider the proof of conservation of mass. By using
$\Psi\equiv 1$ as test function in \eqref{eq:discrete_family}$_1$ we
want to show
\begin{equation}
  \begin{split}
    0 = \ddt \qp{\int_\W \rho_h \d \geovec x}
    &= 
    \int_\W \partial_t \rho_h \d \geovec x 
    =
    -
    \int_\W \div (\rho_h \geovec v_h ) \d \geovec x 
    -
    \int_\E F_1(\rho_h,\geovec v_h, \tau_h, \geovec q_h, 1) \, \d s\\
    &= 
    -
    \int_\E \jump{\rho_h \geovec v_h} \avg{1} \d \geovec x 
    -
    \int_\E F_1(\rho_h,\geovec v_h, \tau_h, \geovec q_h, 1) \, \d s,
  \end{split}
\end{equation}
by Proposition \ref{Pro:trace-jump-avg} with $\geovec p = \rho_h
\geovec v_h$ and $\phi = 1$ and noting $\geovec v_h = \geovec 0$ on
$\partial\W$. Hence the scheme conserves mass if the condition
(\ref{eq:cond_mass}) is true.

Let us now turn to the conservation of energy. Define
\begin{equation}
  E(\rho_h,\geovec v_h, \geovec q_h)
  :=
  \int_\W W(\rho_h) 
  +
  \frac{1}{2}\rho_h \norm{\geovec v_h}^2 
  +
  \frac{\gamma}{2} \norm{\geovec q_h}^2 \d \geovec x.
\end{equation}
Again we want to show
\begin{equation}
  0 
  =
  \ddt {  E(\rho_h,\geovec v_h, \geovec q_h)}.
\end{equation}
Explicitly computing the time derivative
\begin{equation}
  \begin{split}
    \ddt E(\rho_h,\geovec v_h, \geovec q_h)
    &= 
    \int_\W 
    W'(\rho_h) \pdt{\rho}
    +
    \frac{1}{2} 
    \pdt{\rho_h} 
    \norm{\geovec v_h}^2 
    +
    \rho_h \Transpose{\qp{\geovec v_h}}
    \pdt{\geovec v_h}
    +
    \gamma \Transpose{\qp{\geovec q_h}} \pdt{\geovec q_h} \d \geovec x.
  \end{split}
\end{equation}
In view of (\ref{eq:discrete_family})$_4$ and Proposition \ref{Pro:trace-jump-avg} we see
\begin{equation}
  \begin{split}
    \ddt E(\rho_h,\geovec v_h, \geovec q_h)
    &=
    \int_\W
    W'(\rho_h) \pdt{\rho} 
    +
    \frac{1}{2} \pdt{\rho_h} \norm{\geovec v_h}^2 
    +
    \rho_h \Transpose{\qp{\geovec v_h}} \pdt{\geovec v_h}
    +
    \gamma \Transpose{\qp{\nabla \pdt{\rho_h}}} \geovec q_h 
    \d \geovec x 
    \\
    & \qquad \qquad -
    \gamma 
    \int_\E 
    \D_t F_4(\rho_h,\geovec v_h,\tau_h, \geovec q_h, \geovec q_h)
    \d s
    \\
    &=
    \int_\W 
    W'(\rho_h) \pdt{\rho}
    +
    \frac{1}{2} \pdt{\rho_h} \norm{\geovec v_h}^2 
    +
    \rho_h \Transpose{\qp{\geovec v_h}} \pdt{\geovec v_h}
    -
    \gamma \Transpose{\qp{ \pdt{\rho_h}}} \div \geovec q_h 
    \d \geovec x 
    \\
    & \qquad \qquad 
    -\gamma
    \int_\E
    \D_t F_4(\rho_h,\geovec v_h,\tau_h, \geovec q_h, \geovec q_h) 
    -
     \jump{(\pdt{\rho_h}) \geovec q_h}\d s,
  \end{split}
\end{equation}
as $\Transpose{\geovec q_h}\geovec n=0$ on $\partial \W.$
Making use of (\ref{eq:discrete_family})$_2$ and (\ref{eq:discrete_family})$_3$ we see
\begin{equation}
  \begin{split}
    \ddt E(\rho_h,\geovec v_h, \geovec q_h)
    &=
    \int_\W 
    \tau_h \pdt{\rho_h}
    -
    \Transpose{\geovec v_h }\div(\rho_h \geovec v_h \otimes \geovec v_h ) 
    +
    \div(\rho_h \geovec v_h) \norm{\geovec v_h}^2
    - 
    \rho_h \Transpose{\geovec v_h} \nabla \tau_h 
    \\
    &
    \qquad\qquad
    +
    \frac{1}{2} \rho_h \Transpose{\geovec v_h} \nabla \qp{\norm{\geovec v_h}^2}
    \d \geovec x
    \\
    &
    \qquad
    -
    \int_\E 
    \gamma \D_t F_4(\rho_h,\geovec v_h,\tau_h, \geovec q_h, \geovec q_h)
    -
    \gamma \jump{(\pdt{\rho_h}) \geovec q_h}
    \\
    &
    \qquad\qquad\qquad
    -
    F_3(\rho_h,\geovec v_h,\tau_h, \geovec q_h, \pdt{\rho_h}) 
    +
    F_2 (\rho_h,\geovec v_h,\tau_h, \geovec q_h, \geovec v_h)\d s.
  \end{split}
\end{equation}
Now by (\ref{eq:discrete_family})$_1$ and Proposition
\ref{Pro:trace-jump-avg} we see, as $\geovec v_h \in \feszero^d$,
\begin{equation}
  \begin{split}
    \ddt E(\rho_h,\geovec v_h, \geovec q_h)
    &=
    \int_\W 
    - \div(\rho_h \geovec v_h) \tau_h 
    -
    \Transpose{\geovec v_h}\div(\rho_h \geovec v_h \otimes \geovec v_h ) 
    +
    \div(\rho_h \geovec v_h) \norm{\geovec v_h}^2
    \\
    &
    \qquad\qquad
    - 
    \rho_h \geovec v_h \nabla \tau_h
    +
    \frac{1}{2} \rho_h \geovec v_h \nabla \qp{\norm{\geovec v_h}^2} 
    \d \geovec x
    \\
    &
    \qquad    
    -
    \int_\E \gamma \D_t F_4(\rho_h,\geovec v_h,\tau_h, \geovec q_h, \geovec q_h) 
    -
    \gamma \jump{(\pdt{\rho_h}) \geovec q_h}
    \\
    &
    \qquad\qquad\qquad
    -
    F_3(\rho_h,\geovec v_h,\tau_h, \geovec q_h, \pdt{\rho_h})
    \d s
    \\
    &
    \qquad
    -
    \int_\E 
    F_2 (\rho_h,\geovec v_h,\tau_h, \geovec q_h, \geovec v_h)
    + 
    F_1 (\rho_h,\geovec v_h,\tau_h, \geovec q_h, \tau_h)
    \d s
    \\
    &=
    -\int_\E 
    \gamma \D_t F_4(\rho_h,\geovec v_h,\tau_h, \geovec q_h, \geovec q_h) 
    -
    \gamma \jump{(\pdt{\rho_h}) \geovec q_h}
    \\
    &
    \qquad\qquad
    -
    F_3(\rho_h,\geovec v_h,\tau_h, \geovec q_h, \pdt{\rho_h})\d s\\
    & \qquad 
    -
    \int_\E F_2 (\rho_h,\geovec v_h,\tau_h, \geovec q_h, \geovec v_h)
    + 
    F_1 (\rho_h,\geovec v_h,\tau_h, \geovec q_h, \tau_h) 
    \\
    &
    \qquad\qquad\qquad
    + \jump{\rho_h \tau_h \geovec v_h}\d s.
 \end{split}
\end{equation}
Thus, an energy conserving scheme has to satisfy 
\begin{equation}
  \begin{split}
    \label{eq:cond_energy}
    0
    &=
    \int_\E 
    -\gamma  \D_t F_4(\rho_h,\geovec v_h,\tau_h, \geovec q_h, \geovec q_h) 
    +
    \gamma \jump{(\pdt{\rho_h}) \geovec q_h}
    + 
    F_3(\rho_h,\geovec v_h,\tau_h, \geovec q_h, \pdt{\rho_h})\d s
    \\
    & \qquad\qquad 
    -
    \int_\E F_2 (\rho_h,\geovec v_h,\tau_h, \geovec q_h, \geovec v_h)
    +
    F_1 (\rho_h,\geovec v_h,\tau_h, \geovec q_h, \tau_h) 
    +
    \jump{\rho_h \tau_h \geovec v_h}\d s.
 \end{split}
\end{equation}
Note that when \eqref{eq:cond_energy} holds, $F_4$ cannot depend on
$\geovec v_h, \tau_h,\geovec q_h$. Furthermore every summand in $D_t
F_4$ and $F_3$ has to depend on $\pdt{\rho_h}.$ As the trace of
$\pdt{\rho_h}$ is independent of the traces of $\rho_h,\geovec v_h,
\tau_h , \geovec q_h$ the quantities containing $\pdt{\rho_h}$ must
cancel each other, and the ones not containing $\pdt{\rho_h}$ must
cancel each other.
\end{Proof}

\begin{Rem}
  Similarly to the proof of Theorem \ref{The:conservation} one can show
  that  a scheme dissipates energy (even for $\mu=0$), \ie
  \begin{equation}
    \ddt \qp{\int_\W W(\rho_h) 
      + \frac{1}{2}\rho_h \norm{\geovec v_h}^2 
      + \frac{\gamma}{2} \norm{\geovec q_h}^2 \d \geovec x }\leq 0
    \end{equation}
  provided
  \begin{equation}
    \begin{split}
      \int_\E 
      F_1(\rho_h,\geovec v_h, \tau_h, \geovec q_h, \tau_h) 
      +
      F_2(\rho_h,\geovec v_h, \tau_h, \geovec q_h, \geovec v_h) 
      +
      \jump{\rho_h \tau_h \geovec v_h}\d s  &\geq 0  
      \\
      \int_\E
      F_3(\rho_h,\geovec v_h, \tau_h, \geovec q_h, \partial_t \rho_h) 
      -
      \gamma  \D_t F_4(\rho_h,\geovec v_h, \tau_h, \geovec q_h ,\geovec q_h )
      + 
      \gamma \jump{\pdt{ \rho_h} \geovec q_h} \d s &\leq 0.
\end{split}
\end{equation}
This opens the way for the construction of schemes which dissipate a
small amount of energy, which might serve as a stabilising mechanism
for the scheme, \eg in case forward time stepping is considered.
\end{Rem}

\begin{Example}
  For $\alpha,\beta>0$ the choice of fluxes
  \begin{equation}
    \begin{split}
      F_1(\rho_h, \geovec v_h, \tau_h, \geovec q_h , \Psi)
      &=
      -
      \jump{\rho_h \geovec v_h} 
      \avg{\Psi}
      +
      \alpha 
      \Transpose{\jump{\tau_h}}
      \jump{\Psi}
      \\
      F_2(\rho_h, \geovec v_h, \tau_h, \geovec q_h ,\geovec \Chi)
      &=
      -
      \Transpose{\jump{\tau_h}} 
      \avg{\rho_h \geovec \Chi}
      +
      \beta
      \jump{\geovec  v_h} 
      \jump{\geovec \Chi}
      \\
      F_3(\rho_h, \geovec v_h, \tau_h, \geovec q_h , \Xi)
      &=
      -
      \gamma
      \jump{\geovec q_h}
      \avg{\Xi}
      \\
      F_4(\rho_h, \geovec v_h, \tau_h, \geovec q_h ,\geovec \Zeta)
      &=
      \Transpose{\jump{\rho_h}} 
      \avg{\geovec \Zeta},
    \end{split}
  \end{equation}
  in \eqref{eq:discrete_family} yields a scheme which is consistent,
  conserves mass and dissipates energy.
\end{Example}

\begin{Proof}[of Corollary \ref{cor:energy-dissipation}]
 The proof of conservation of mass is exactly the same as in the proof
 of Theorem \ref{The:conservation} because
 \eqref{eq:discrete_family}$_1$ does not depend on $\mu.$ For the
 dissipation of energy the only difference to the proof of Theorem
 \ref{The:conservation} is that, when \eqref{eq:discrete_family}$_2$
 is tested with $\geovec \Chi=\geovec v_h$ an additional summand $\mu
 B_h(\geovec v_h,\geovec v_h)$ is created and this term is not altered
 by the subsequent calculations.
\end{Proof}


We will now show that it is very restrictive for schemes given by
\eqref{eq:discrete_family} to be consistent, conserve 
momentum and energy.

\begin{Pro}
  \label{the:mom-cons-cond}
  Let $\geovec e_i$ be the $i$--th coordinate vector of $\reals^d$. A
  generic scheme of the form \eqref{eq:discrete_family} is momentum
  conservative, \ie satisfies
  \begin{equation}
    \ddt \int_\W \rho_h \geovec v_h \d \geovec x = \geovec 0
  \end{equation}
  if and only if
  \begin{equation}
    \label{eq:momentum_change}
    \begin{split}
      0 &=
      -
      \int_\E
      F_1(\rho_h,\geovec v_h,\tau_h,\geovec q_h ,\Transpose{\geovec e_i} \geovec v_h )
    +
    F_2(\rho_h,\geovec v_h,\tau_h,\geovec q_h ,{\geovec e_i}) 
    \\
    &
    \qquad\qquad\qquad
    +
    F_3(\rho_h,\geovec v_h,\tau_h,\geovec q_h ,\partial_{x_i}\rho_h) 
    +
    F_4(\rho_h,\geovec v_h,\tau_h,\geovec q_h ,\nabla \partial_{x_i}\rho_h) 
    \d s
    \\
    &\qquad 
    +
    \int_{\E \cup \partial \W}
    - 
    \jump{\rho_h \tau_h \geovec e_i}
    +
    \jump{\qp{W(\rho_h) + \frac{1}{2} \rho_h \norm{\geovec v_h}^2} \geovec e_i}
    -
    \jump{\rho_h \Transpose{\geovec e_i} \qp{\geovec
  v_h\otimes \geovec v_h 
      }}
    \\
    &
    \qquad\qquad\qquad
    -
    \gamma 
    \jump{\partial_{x_i} \rho_h \geovec q_h} 
    +
    \gamma
    \jump{\rho_h \nabla \partial_{x_i} \rho_h} 
    \\
    &
    \qquad\qquad\qquad
    -
      \gamma
      \jump{
        \rho_h \Delta \rho_h \geovec e_i 
        +
      \frac{1}{2} 
      \norm{\nabla \rho_h}^2 
      \geovec e_i
      -
      \partial_{x_i} \rho_h 
      \nabla \rho_h
    }
      \d s.
  \end{split}
\end{equation}

\end{Pro}
\begin{Proof}
  We start with a rather general calculation of the change of
  momentum in direction $\geovec e_i$ for $i=1,\dots,d$ for the
  generic scheme (\ref{eq:discrete_family}). We define
  \begin{equation}
    M_i(\rho_h,\geovec v_h)
    :=
    \int_\W 
    \rho_h 
    \Transpose{\geovec e_i}
    \geovec v_h
    \d \geovec x
  \end{equation}
  to be the discrete momentum in direction $\geovec e_i$.
  
  Then the rate of change of the discrete momentum is
  \begin{equation}
    \ddt 
    M_i(\rho_h, \geovec v_h)
    =
    \int_\W 
    \pdt{\rho_h} 
    \Transpose{\geovec e_i}
    \geovec v_h
    +
    \Transpose{\geovec e_i}
    \rho_h 
    \pdt{\geovec v_h}
    \d \geovec x.
  \end{equation}
  Now making use of (\ref{eq:discrete_family})$_1$ and
  (\ref{eq:discrete_family})$_2$ use $\Psi = \Transpose{\geovec e_i}\geovec v_h$ and $\geovec \Chi =
  \geovec{e}_i$,
  \begin{equation}
    \begin{split}
      \ddt 
      M_i(\rho_h, \geovec v_h)
      &=
      \int_\W 
      -
      \div(\rho_h \geovec v_h)
      \Transpose{\geovec e_i}
      \geovec v_h 
      -
      \div\qp{
        \rho_h  \Transpose{\geovec e_i}
        \qp{\geovec v_h \otimes \geovec v_h}
      } 
      + 
      \div\qp{
        \rho_h \geovec v_h
      }
      \Transpose{\geovec e_i}
      \geovec v_h
      \\
      &
      \qquad\qquad
      -
      \rho_h 
      \Transpose{\geovec e_i}
      \nabla \tau_h
      +
      \frac{1}{2}\rho_h 
      \Transpose{\geovec e_i}
      \nabla \qp{\norm{\geovec v_h}^2}
      \d \geovec x 
      \\
      &
      \qquad
      -
      \int_\E 
      F_1(\rho_h,\geovec v_h,\tau_h,\geovec q_h ,\Transpose{\geovec e_i} \geovec v_h)
      +
      F_2(\rho_h,\geovec v_h,\tau_h,\geovec q_h ,{\geovec e_i})
      \d s
      .
    \end{split}
  \end{equation}
  Using (\ref{eq:discrete_family})$_3$ together with Proposition
  \ref{Pro:trace-jump-avg},
\begin{equation}
  \begin{split}
    \ddt 
    M_i(\rho_h, \geovec v_h)
    &=
    \int_\W 
    -
    \div\qp{
      \rho_h  
      \Transpose{\geovec e_i}
      \qp{\geovec v_h \otimes 
        \geovec v_h}
    }
    +
    \div\qp{\rho_h \geovec e_i}
    \tau_h
    +
    \frac{1}{2}
    \rho_h 
    \Transpose{\geovec e_i} 
    \nabla \qp{\norm{\geovec v_h}^2}
    \d \geovec x 
    \\
    &
    \qquad 
    -
    \int_\E 
    F_1(\rho_h,\geovec v_h,\tau_h,\geovec q_h ,\Transpose{\geovec e_i} \geovec v_h ) 
    +
    F_2(\rho_h,\geovec v_h,\tau_h,\geovec q_h ,{\geovec e_i})
    \d s
    \\
    &
    \qquad\qquad\qquad
    - 
    \int_{\E \cup \partial \W}
    \jump{\rho_h \tau_h \geovec e_i} 
    \d s
    \\
    &=
    \int_\W 
    -
    \div\qp{
      \rho_h 
      \Transpose{\geovec e_i} 
      \qp{
        \geovec v_h \otimes \geovec v_h
      }
    }
    +
    \div\qp{
      W(\rho_h)
      \geovec e_i
    }
    \\
    &
    \qquad\qquad
    -
    \gamma 
    \div\qp{\geovec q_h}
    \div\qp{\rho_h \geovec e_i} 
    + 
    \frac{1}{2} 
    \div\qp{
      \rho_h
      \norm{\geovec v_h}^2 \geovec e_i
    }
    \d \geovec x
    \\
    &
    \qquad
    -
    \int_\E
    F_1(\rho_h,\geovec v_h,\tau_h,\geovec q_h ,\Transpose{\geovec e_i} \geovec v_h ) 
    +
    F_2(\rho_h,\geovec v_h,\tau_h,\geovec q_h ,{\geovec e_i})
    \\
    &
    \qquad\qquad\qquad
    +
    F_3(\rho_h,\geovec v_h,\tau_h,\geovec q_h ,\div(\rho_h \geovec e_i)) 
    \d s
    -
    \int_{\E \cup \partial \W}
    \jump{\rho_h \tau_h \geovec e_i} 
    \d s.
  \end{split}
\end{equation}
Again, in view of Proposition \ref{Pro:trace-jump-avg} we may
integrate by parts elementwise and it follows
\begin{equation}
  \begin{split}
    \ddt 
    M_i(\rho_h, \geovec v_h)
    &=
    \int_\W 
    \gamma 
    \Transpose{ \geovec q_h }
    \nabla \partial_{x_i} \rho_h 
    \d \geovec x
    -
    \int_\E
    F_1(\rho_h,\geovec v_h,\tau_h,\geovec q_h ,\Transpose{\geovec e_i} \geovec v_h ) 
    +
    F_2(\rho_h,\geovec v_h,\tau_h,\geovec q_h ,{\geovec e_i})
    \\
    &
    \qquad\qquad 
    +
    F_3(\rho_h,\geovec v_h,\tau_h,\geovec q_h ,\div(\rho_h \geovec e_i)) 
    \d s
    \\
    &
    \qquad 
    +
    \int_{\E \cup \partial\W}
    -
    \jump{\rho_h \tau_h \geovec e_i}
    +
    \jump{\qp{W(\rho_h) + \frac{1}{2} \rho_h \norm{\geovec v_h}^2}\geovec e_i}
    \\
    &
    \qquad\qquad\qquad 
    -
    \jump{\rho_h \Transpose{\geovec e_i} \qp{\geovec v_h \otimes \geovec v_h }}
    - 
    \gamma \jump{\partial_{x_i} \rho_h  \geovec q_h}
    \d s.
  \end{split}
\end{equation}
By (\ref{eq:discrete_family})$_4$ and another application of
Proposition \ref{Pro:trace-jump-avg} it holds that
\begin{equation}
  \label{eq:mom-cons-proof-1}
  \begin{split}
    \ddt 
    M_i(\rho_h, \geovec v_h)
    &=
    \int_\W 
    -
    \gamma 
    \rho_h 
    \Delta \partial_{x_i} 
    \rho_h 
    \d \geovec x
    -
    \int_\E 
    F_1(\rho_h,\geovec v_h,\tau_h,\geovec q_h ,\Transpose{\geovec e_i} \geovec v_h )
    +
    F_2(\rho_h,\geovec v_h,\tau_h,\geovec q_h ,{\geovec e_i})
    \\
    &
    \qquad\qquad
    +
    F_3(\rho_h,\geovec v_h,\tau_h,\geovec q_h ,\partial_{x_i}\rho_h) 
    +
    F_4(\rho_h,\geovec v_h,\tau_h,\geovec q_h ,\nabla \partial_{x_i}\rho_h) 
    \d s
    \\
    &
    \qquad 
    +
    \int_{\E \cup \partial \W} 
    -
    \jump{\rho_h \tau_h \geovec e_i} 
    +
    \jump{\qp{W(\rho_h) + \frac{1}{2} \rho_h \norm{\geovec v_h}^2}\geovec e_i}
    -
    \jump{\rho_h \Transpose{\geovec e_i} \qp{\geovec v_h \otimes \geovec v_h }}
    \\
    &
    \qquad\qquad\qquad
    - 
    \gamma
    \jump{\partial_{x_i} \rho_h  \geovec q_h}
    +
    \gamma 
    \jump{\rho_h \nabla \partial_{x_i} \rho_h}
    \d s.
  \end{split}
\end{equation}
We may write the first integral appearing on the right hand side of
(\ref{eq:mom-cons-proof-1}) in the following way:
\begin{equation}
  \begin{split}
    \ddt 
    M_i(\rho_h, \geovec v_h)
    &=
    \int_\W 
    -
    \gamma 
    \div\qp{
      \rho_h 
      \Delta \rho_h 
      \geovec e_i 
      +
      \frac{1}{2} 
      \norm{\nabla \rho_h}^2
      \geovec e_i 
      -
      \partial_{x_i}
      \rho_h 
      \nabla \rho_h
      }
    \d \geovec x
    \\
    &\qquad 
    -
    \int_\E
    F_1(\rho_h,\geovec v_h,\tau_h,\geovec q_h ,\Transpose{\geovec e_i} \geovec v_h)
    +
    F_2(\rho_h,\geovec v_h,\tau_h,\geovec q_h ,{\geovec e_i})
    \\
    &
    \qquad\qquad\qquad
    +
    F_3(\rho_h,\geovec v_h,\tau_h,\geovec q_h ,\partial_{x_i}\rho_h) 
    +
    F_4(\rho_h,\geovec v_h,\tau_h,\geovec q_h ,\nabla\partial_{x_i}\rho_h)
    \d s
    \\
    &
    \qquad 
    +
    \int_{\E \cup \partial \W} 
    -
    \jump{\rho_h \tau_h \geovec e_i}
    +
    \jump{\qp{W(\rho_h) 
      +
      \frac{1}{2} 
      \rho_h
      \norm{\geovec v_h}^2
    }\geovec e_i}
    -
    \jump{\rho_h \Transpose{\geovec e_i} \qp{\geovec v_h\otimes\geovec v_h }}
    \\
    &
    \qquad\qquad\qquad
    -
    \gamma 
    \jump{\partial_{x_i} \rho_h \geovec q_h}
    +
    \gamma 
    \jump{\rho_h \nabla \partial_{x_i} \rho_h}
    \d s.
  \end{split}
\end{equation}
Hence we may apply Proposition \ref{Pro:trace-jump-avg} one more time
yielding the desired result.

\end{Proof}

\begin{Rem}
  Since the right hand side of (\ref{eq:momentum_change}) depends
  nonlinearly on $\rho_h$, it is not expected that the jump terms
  involving $\nabla \rho_h$ and $\geovec q_h$ will necessarily cancel
  with the other terms appearing. Only in the case $F_4\equiv 0$ will
  the $\nabla \rho_h$ and $\geovec q_h$ terms cancel with each
  other. In this case (\ref{eq:energy_cons})$_2$ gives us the
  condition that
  \begin{equation}
    F_3(\rho_h,\geovec v_h,\tau_h,\geovec q_h ,\Xi) = -\gamma \jump{\geovec q_h \Xi}
  \end{equation}
  which excludes consistency. In the event that $F_4 \not\equiv 0$, the
  terms involving $\nabla \rho_h$ in (\ref{eq:momentum_change}) would
  be required to cancel independently, yielding a condition on $F_3$
  which again excludes consistency.  
\end{Rem}


\section{Development of energy consistent numerical methods -- the temporal discrete case}
\label{sec:num-methods-temporal}

For the readers convenience we will present argument for designing the
temporally discrete {scheme} in the spatially continuous setting. To
obtain a fully discrete version the spatial and temporal
discretisations have to be combined which is straightforward as
presented in \eqref{eq:fully_discrete} and Theorem
\ref{the:fully_discrete}.

We subdivide the time interval $[0,T]$ into a partition of $N$
consecutive adjacent subintervals whose endpoints are denoted
$t_0=0<t_1<\ldots<t_{N}=T$.  The $n$-th timestep is defined as ${k_n :=
t_{n+1} - t_{n}}$.  We will consistently use the shorthand
$F^n(\cdot):=F(\cdot,t_n)$ for a generic time function $F$. We also
denote $\nplush{F} := \frac{1}{2}\qp{F^n + F^{n+1}}$.

\begin{The}
  \label{The:time-discrete-cons}
  Given initial conditions $\rho^0, \geovec v^0, \tau^0 \AND \geovec
  q^0$ the temporal semi discrete scheme is: For $n\in\naturals$, find
  $\rho^{n+1}, \geovec v^{n+1} , \tau^{n+1} \AND \geovec q^{n+1}$ such
  that $\geovec v^{n+1} = \geovec 0 \AND \Transpose{\qp{\geovec
      q^{n+1}}}\geovec n = 0$ on $\partial \W$ and
  \begin{equation}
    \label{eq:time_discrete}
    \begin{split}
      0
      &=
      \frac{\rho^{n+1} - \rho^n}{k_n} 
      +
      \div \qp{\nplush{\rho} \nplush{ \geovec v}}
      \\
      0 
      &=
      \nplush{\rho} 
      \qp{
        \frac{\geovec v^{n+1}
          -
          \geovec v^n}{k_n}
      }
      + 
      \div
      \qp{
        \nplush{\rho}
        \nplush{\geovec v}
        \otimes
        \nplush{\geovec v}
      }
      - 
      \div 
      \qp{
        \nplush{\rho} 
        \nplush{\geovec v} 
      }
      \nplush{\geovec v}
      \\
      &
      \qquad \qquad
      +
      \nplush{\rho_h} 
      \nabla \nplush{\tau}
      -
      \frac{1}{2} 
      \nplush{\rho}
      \nabla 
      \qp{ \norm{\nplush{\geovec v}}^2} - \mu \Delta \nplush{\geovec v}
      \\
      0
      &=
      \nplush{\tau} 
      -
      \frac{
        W(\rho^{n+1})
        - 
        W(\rho^n)}
           {\rho^{n+1} - \rho^n} 
           +
           \gamma 
           \div\qp{
             \nplush{\geovec q}}
           -
           \frac{1}{4}
           \qp{
             \norm{\geovec v^{n+1}}^2
             +
             \norm{\geovec v^{n}}^2}
           \\
           0 
           &=
           \geovec  q^{n+1} 
           -
           \nabla \rho^{n+1}.
    \end{split}
  \end{equation}

  This scheme satisfies the following energy dissipation property for all $0 \leq n \leq N$
  \begin{equation}
    \begin{split}
      \int_\Omega
    W(\rho^n) 
    +
    \frac{1}{2}
    \rho^n 
    \norm{\geovec v^n}^2 
    +
    \frac{\gamma}{2} 
    \norm{\geovec q^n}^2
    \d \geovec x 
    &=
    \int_\Omega 
    W(\rho^0) 
    +
    \frac{1}{2} \rho^0 
    \norm{\geovec v^0}^2 
    +
    \frac{\gamma}{2} 
    \norm{\geovec q^0}^2 
    \d \geovec x 
    \\
    &
    \qquad\qquad
    -
    \mu \sum_{j=0}^{n-1}  k_j \int_\W \norm{\D \jplush{\geovec v}}^2 \d \geovec x.
    \end{split}
  \end{equation}
\end{The}
\begin{Proof}
  We proceed by multiplying \eqref{eq:time_discrete}$_1$ by
  $\nplush{\tau}$ and \eqref{eq:time_discrete}$_2$ by $\nplush{\geovec
    v}$, integrate over the domain $\W$ and take the sum. We obtain
  \begin{equation}
    \label{eq:time_discrete_all}
    0
    =
    \int_\W  \cI_1 + \cI_2 + \cI_3  + \cI_4
    \d \geovec x
  \end{equation}
  with
  \begin{align}
    \cI_1 
    &:=
    \frac{\rho^{n+1} - \rho^n}{k_n}
    \qp{\frac{W(\rho^{n+1}) - W(\rho^n)}{\rho^{n+1} - \rho^n} 
      -
      \gamma \div\qp{\nplush{\geovec q}}
      +
      \frac{1}{4}\qp{\norm{\geovec v^{n+1}}^2 + \norm{\geovec v^{n}}^2}   
    }
    \\
    \nonumber
    &
    \qquad\qquad\qquad
    +
    \nplush{\rho}
    \Transpose{\qp{\nplush{\geovec v}}}
    \qp{
      \frac{\geovec v^{n+1} - \geovec v^n}{k_n}}
    \\
    \cI_2 
    &:=
    \div{\qp{
      \nplush{\rho} 
      \nplush{ \geovec v}}
    }
    \nplush{\tau}
    +
    \nplush{\rho} 
    \Transpose{\qp{\nplush{\geovec v}}}
    \nabla \nplush{\tau}
    \\
    \cI_3 
    &:=
    \Transpose{\qp{\nplush{\geovec v}}}\div \qp{
      \nplush{\rho} 
      \nplush{\geovec v}
      \otimes 
      \nplush{\geovec v}
    }    
    - 
    \div \qp{
      \nplush{\rho} 
      \nplush{\geovec v} 
    }
    \norm{\nplush{\geovec v}}^2
    \\
    & 
    \nonumber\qquad\qquad\qquad
    - 
    \frac{1}{2}
    \nplush{\rho} 
    \Transpose{\qp{\nplush{\geovec v}}}
    \nabla \qp{ \norm{\nplush{\geovec v}}^2} \\
\cI_4 &:= - \mu \int_\W \Transpose{\qp{\nplush{\geovec v}}}\qp{\Delta \nplush{\geovec v}} \d \geovec x.
\end{align}
One may readily check that
\begin{equation}
  \label{eq:time_discrete_I1}
  \begin{split}
     k_n \int_\W 
   \cI_1  
    \d \geovec x
    &=
    \int_\W 
    W(\rho^{n+1}) 
    + 
    \frac{1}{2} 
    \rho^{n+1}
    \norm{\geovec v^{n+1}}^2
    + 
    \frac{\gamma}{2} 
    \norm{\geovec q^{n+1}}^2
    \d \geovec x
    \\
    &
    \qquad
    -  
    \int_\W 
    W(\rho^n)
    +
    \frac{1}{2} 
    \rho^n 
    \norm{\geovec v^n}^2
    +
    \frac{\gamma}{2} 
    \norm{\geovec q^n}^2 
    \d \geovec x
    \\
    &  
    \qquad
    -
    \gamma 
    \int_{\partial \W} 
    (\rho^{n+1} -\rho^n)
    \Transpose{(\geovec q^{n+1} + \geovec q^n)}
    \geovec n 
    \d s
    \\
    &= 
    \int_\W 
    W(\rho^{n+1})
    + 
    \frac{1}{2} 
    \rho^{n+1}
    \norm{\geovec v^{n+1}}^2
    +
    \frac{\gamma}{2}
    \norm{\geovec q^{n+1}}^2 
    \d \geovec x
    \\
    &
    \qquad
    -
    \int_\W  W(\rho^n)
    +
    \frac{1}{2} 
    \rho^n 
    \norm{\geovec v^n}^2 
    +
    \frac{\gamma}{2}
    \norm{\geovec q^n}^2 
    \d \geovec x
 \end{split}
\end{equation}
Moreover 
\begin{equation}
  \label{eq:time_discrete_I2}
  \int_\W  
  \cI_2 
  \d \geovec x 
  =
  \div \qp{
    \nplush{\rho}
    \nplush{ \geovec v}
    \nplush{\tau} 
  }
  \d \geovec x
  =
  \int_{\partial \W} 
  \nplush{\rho}
  \nplush{\tau}
  \Transpose{\qp{\nplush{ \geovec v}} }
  \geovec n
  \d s
  =
  0
  .
\end{equation}
Furthermore we see that $\cI_3$ satisfies
\begin{equation}
  \label{eq:time_discrete_I3}
  \begin{split}
    \cI_3 
    &=
    \sum_{i,j=1}^d 
    \partial_{x_i} 
    \qp{
      \nplush{\rho}
      \nplush{v_i} 
      \nplush{v_j}}
    \nplush{v_j}
    -
    \partial_{x_i}
    \qp{
      \nplush{\rho} 
      \nplush{v_i}}
    \qp
       {
         \nplush{v_j}
       }^2
       \\
       &
       \qquad\qquad
       -
     \frac{1}{2}
     \nplush{\rho} 
     \nplush{ v_i}
     \partial_{x_i}
     \qp{\qp{
         \nplush{v_j}}^2}
     =
     0.
  \end{split}
\end{equation}
Finally we find for $\cI_4$
\begin{equation}
 \begin{split}
  \int_\W 
    \cI_4 &= \mu \int_\W \norm{\D \nplush{\geovec v}}^2 \d \geovec x 
- \mu \int_{\partial \W} \Transpose{\qp{\nplush{\geovec v}}} \qp{\D \nplush{\geovec v}} \geovec n \d s\\
 &= \mu \int_\W \norm{\D \nplush{\geovec v}}^2 \d \geovec x.
 \end{split}
\end{equation}

Inserting \eqref{eq:time_discrete_I1}--\eqref{eq:time_discrete_I3}
into \eqref{eq:time_discrete_all} yields
\begin{equation}
  \begin{split}
    0
    &=
    \int_\W 
    W(\rho^{n+1}) 
    +
    \frac{1}{2} \rho^{n+1}
    \norm{\geovec v^{n+1}}^2 
    +
    \frac{\gamma}{2} \norm{\geovec q^{n+1}}^2 
    \d \geovec x
   \\& \qquad \qquad - 
    \int_\W
    W(\rho^n) 
    +
    \frac{1}{2}
    \rho^n
    \norm{\geovec v^n}^2 
    +
    \frac{\gamma}{2} \norm{\geovec q^n}^2  
    \d \geovec x    
+ \mu k_n \int_\W \norm{\D \nplush{\geovec v}}^2 \d \geovec x
,
\end{split}
\end{equation}
concluding the proof.
\end{Proof}

\section{Development of consistent numerical methods -- the fully discrete case}
\label{sec:num-methods-fully}
In this section we combine our spatial and temporal discretisations to
provide a fully discrete numerical method for the \EK and \NSK
systems.

Let $\ltwoproj{}:\sobh{1}(\T{}) \to \fes,
\ltwoprojzero{}:\sobh1_0(\T{})^d\to{\feszero}^d \AND
\ltwoprojn{}:\hon(\T{})\to\fesn$ be the $\leb{2}$ projection operators
into $\fes, {\feszero}^d \AND \fesn$ respectively. We combine the
arguments given in \S\ref{sec:numerical-methods} and
\S\ref{sec:num-methods-temporal} to obtain a fully discrete scheme which,
given 
\begin{equation}
  \rho_h^0 := \ltwoproj{\rho^0},
  \quad
  \geovec v_h^0 :=
  \ltwoprojzero{\geovec v^0}, 
  \quad 
  \tau_h^0 := \ltwoproj{\tau^0} 
  \AND\quad \geovec
  q_h^0 := \ltwoprojn{\geovec q^0}
  \end{equation}
requires us to find a sequence of
functions $\rho_h^{n+1} ,\tau_h^{n+1} \in \fes, \geovec
v_h^{n+1} \in {\feszero}^d \AND  \geovec q_h^{n+1}\in\fesn$ such that
\begin{equation}
  \begin{split}
    \label{eq:fully_discrete}
    0 
    &=
    \int_\W \qp{\frac{\rho_h^{n+1}-\rho_h^n}{k_n} + {\div\qp{\rho_h^{n+\frac{1}{2}} \geovec v_h^{n+\frac{1}{2}} }}} {\Psi} \d \geovec x
    \\
    &
    \qquad
    + 
    \int_\E
    F_1\qp{\rho_h^{n+\frac{1}{2}} , \geovec v_h^{n+\frac{1}{2}} , \tau_h^{n+\frac{1}{2}} , \geovec q_h^{n+\frac{1}{2}} , \Psi} \d s
    \Foreach \Psi \in \fes
    \\
    0
    &=
    \int_\W
    \Transpose{\qp{\qp{\rho_h^{n+\frac{1}{2}}  \frac{\geovec v_h^{n+1}-\geovec v_h^n}{k_n} }
      +
      \div\qp{\rho_h^{n+\frac{1}{2}}  \geovec v_h^{n+\frac{1}{2}} \otimes \geovec v_h^{n+\frac{1}{2}} }}}\geovec \Chi \d \geovec x\\
    & 
    \qquad 
    +  \int_\W 
      \Transpose{\qp{
      -
      \div\qp{\rho_h^{n+\frac{1}{2}}  \geovec v_h^{n+\frac{1}{2}} } \geovec v_h^{n+\frac{1}{2}} 
      +
      \rho_h^{n+\frac{1}{2}} \nabla \tau_h^{n+\frac{1}{2}} 
      -
      \frac{1}{2}
      \rho_h^{n+\frac{1}{2}}  \nabla \norm{\geovec v_h^{n+\frac{1}{2}} }^2 
    }} \geovec \Chi
    \d \geovec x
    \\
    &
    \qquad 
    + 
    \int_\E 
    F_2\qp{\rho_h^{n+\frac{1}{2}}, \geovec v_h^{n+\frac{1}{2}}, \tau_h^{n+\frac{1}{2}}, \geovec q_h^{n+\frac{1}{2}}, \geovec \Chi} \d s
    + \mu B_h (\geovec v_h^{n+\frac{1}{2}},\geovec \Chi)
    \Foreach \geovec \Chi \in \feszero^d
    \\
    0 
    &=
    \int_\W
    \qp{  \tau_h^{n+\frac{1}{2}} - \frac{W(\rho_h^{n+1})-W(\rho_h^n)}{\rho_h^{n+1}-\rho_h^n} 
+ \gamma \div\qp{\geovec q_h^{n+\frac{1}{2}}} - \frac{1}{4}\qp{\norm{\geovec v_h^{n+1}}^2 + \norm{\geovec v_h^n}^2 } } \Xi \d \geovec x
    \\
    &
    \qquad 
    +
    \int_\E
    F_3\qp{\rho_h^{n+\frac{1}{2}}, \geovec v_h^{n+\frac{1}{2}}, \tau_h^{n+\frac{1}{2}}, \geovec q_h^{n+\frac{1}{2}}, \Xi} \d s 
    \Foreach \Xi \in \fes
    \\
    0 
    &=
    \int_\W \Transpose{\qp{\geovec q_h^{n+1} - \nabla \rho_h^{n+1}}} \geovec \Zeta \d \geovec x
    +
    \int_\E  F_4\qp{\rho_h^{n+1}, \geovec v_h^{n+1}, \tau_h^{n+1}, \geovec q_h^{n+1}, \geovec \Zeta} \d s
    \Foreach \geovec \Zeta \in \fesn.
  \end{split}
\end{equation}

\begin{The}
  \label{the:fully_discrete}
  Under the assumptions on the fluxes (\ref{eq:cond_mass}),
  (\ref{eq:energy_cons}) given in Theorem \ref{The:conservation}, the
  solution of the scheme \eqref{eq:fully_discrete} conserves mass, \ie
\[ \int_\W \rho_h^n \d \geovec x = \int_\W \rho_h^0 \d \geovec x  \qquad \text{ for } \ 0 \leq  n \leq N\]
and satisfies the energy dissipation equality
\begin{multline}
  \int_\W
  W(\rho_h^{n+1}) 
  +
  \frac{1}{2}\rho_h^{n+1} \norm{\geovec v_h^{n+1}}^2 
  +
  \frac{\gamma}{2} \norm{\geovec q_h^{n+1}}^2 \d \geovec x
  \\
  -
  \int_\W
  W(\rho_h^{n}) 
  +
  \frac{1}{2}\rho_h^{n} 
  \norm{\geovec v_h^{n}}^2 
  +
  \frac{\gamma}{2}
  \norm{\geovec q_h^{n}}^2 \d \geovec x
  =
  -\mu k_n B_h(\nplush{\geovec v_h},\nplush{\geovec v_h} ).
\end{multline}
\end{The}
\begin{Proof}
  The proof is merely combining the results of Theorem
  \ref{The:conservation}, Corollary \ref{cor:energy-dissipation} and
  Theorem \ref{The:time-discrete-cons}.
\end{Proof}


\section{Numerical experiments}
\label{sec:numerics}

In this section we conduct a series of numerical experiments aimed at
testing the robustness of the method. There are four experiments which
investigate the behaviour of the discrete energy for the \EK
(\S\ref{sec:energy-EK}) and the \NSK systems (\S\ref{sec:energy-NSK}),
benchmarking the algorithm against a travelling wave solution of the
\EK system (\S\ref{sec:benchmark}), observing that there are no
parasitic currents in long time simulations. Moreover, we conduct some
simulations for $d=2$ (\S\ref{sec:2d-sims}).

In each of these experiments we consider the fully discrete scheme
(\ref{eq:fully_discrete}) with the numerical fluxes given in Corollary
\ref{cor:consistency-and-conservation}.

\subsection{Implementation issues}

The numerical experiments were conducted using the \dolfin interface
for \fenics \cite{LoggWells:2010}. The graphics were generated using
\gnuplot and \paraview.

In each of the numerical experiments we fix $W$ to be the following
quartic double well potential
\begin{equation}
  \label{eq:double-well}
  W(\rho) = \frac{1}{4}\qp{\rho-1}^2\qp{\rho-2}^2
\end{equation}
with minima at $\rho = 1$ and $\rho = 2$.

\begin{Rem}[the quotient of the double well]
  In the computational implementation we did not use the difference
  quotient $\tfrac{W(\rho^{n+1})-W(\rho^n)}{\rho^{n+1} - \rho^n}$
  appearing in (\ref{eq:time_discrete}) as it is ill-defined for
  $\rho^{n+1} = \rho^n$ and badly conditioned when $|\rho^{n+1} -
  \rho^n|$ is small. Instead we use a sufficiently high order
  approximation of this term. For (\ref{eq:double-well}) we use the
  following Taylor expansion representation
  \begin{equation}
    \frac{W(\rho^{n+1})-W(\rho^n)}{\rho^{n+1} - \rho^n}
    =
    W'(\nplush{\rho})
    +
    \tfrac{1}{24}W'''(\nplush{\rho})\qp{\rho^{n+1} - \rho^{n}}^2
  \end{equation}
  which is exact. We note that when $W$ is not polynomial a
  sufficiently high order truncation of the Taylor expansion can be
  achieved such that the change in energy is of high order with
  respect to the timestep. This allows the
  construction of a method with arbitrarily small deviations of the
  energy with respect to the timestep.
\end{Rem}

In each of the subsequent numerical experiments we assemble the
discrete system (\ref{eq:fully_discrete}) as a nonlinear system of equations. The solution
to the nonlinear system was approximated by a Newton solver with a
tolerance set to $10^{{-10}}$. On each Newton step the linear system
of equations was approximated using a stabilised conjugate gradient
solver with an incomplete LU preconditioner also set to a tolerance of
$10^{{-10}}$.

\subsection{Test 1 -- conservativity for the \EK system}
\label{sec:energy-EK}

In this case we take $\mu = 0$. We are then studying the
conservativity property of the numerical method proposed for the \EK
system in Corollary \ref{cor:consistency-and-conservation}. We take
$\W = [0,1]$ and consider an initial condition given by a step
function
\begin{equation}
  \label{eq:step-ics}
  \rho_0(x) = 
  \begin{cases}
    1.1 \text{ if } x \leq 0.5
    \\
    1.9 \text{ otherwise },
  \end{cases}
\qquad 
v_0 \equiv 0.
\end{equation}

\begin{figure}[ht]
  \caption[]
          {
            \label{fig:1d-conservative-energy-plot} 
            \ref{sec:energy-EK} Test 1 -- Numerical experiment showing
            the conservation of mass and energy for the numerical
            method proposed in Corollary
            \ref{cor:consistency-and-conservation} for the \EK system
            (\ie $\mu = 0$).  Due to the energy conservativity the \EK
            simulation will never achieve a steady state, the
            oscillations will continue to propogate.}
          \begin{center}
            \subfigure[][Initial condition, $t = 0$]{
              \includegraphics[scale=\figscale, width=0.47\figwidth]
                              {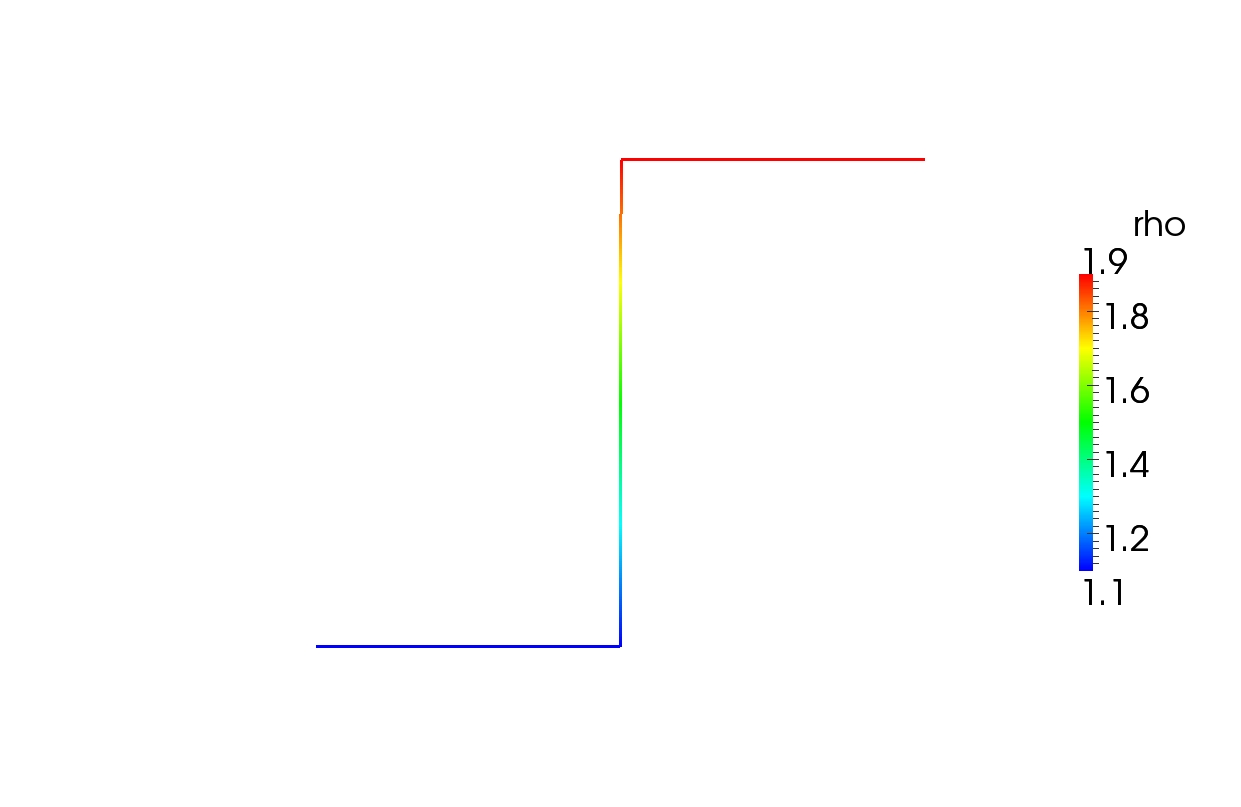} }
            \subfigure[][$t = 0.01$]{ \includegraphics[scale=\figscale,
                width=0.47\figwidth] {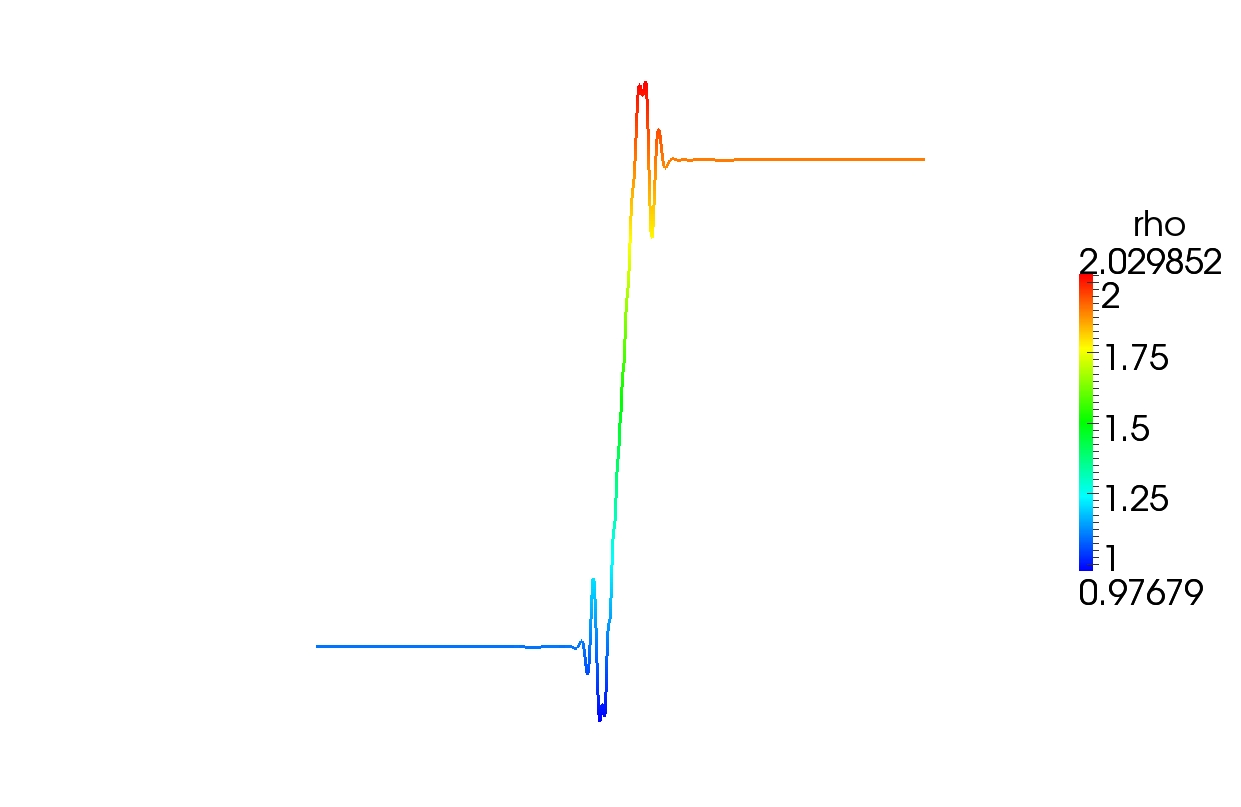}
            } \subfigure[][$t = 0.05$]{ \includegraphics[scale=\figscale,
                width=0.47\figwidth] {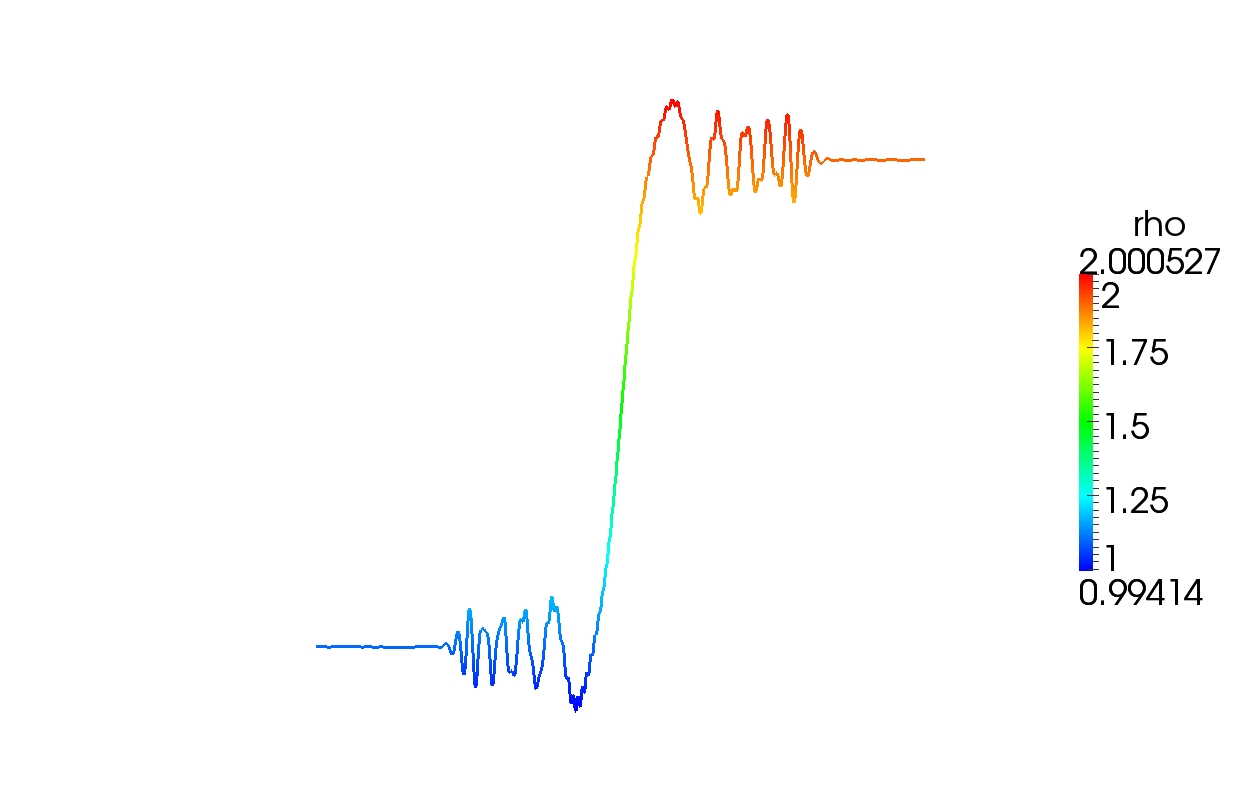}
            } \subfigure[][$t = 0.1$]{ \includegraphics[scale=\figscale,
                width=0.47\figwidth] {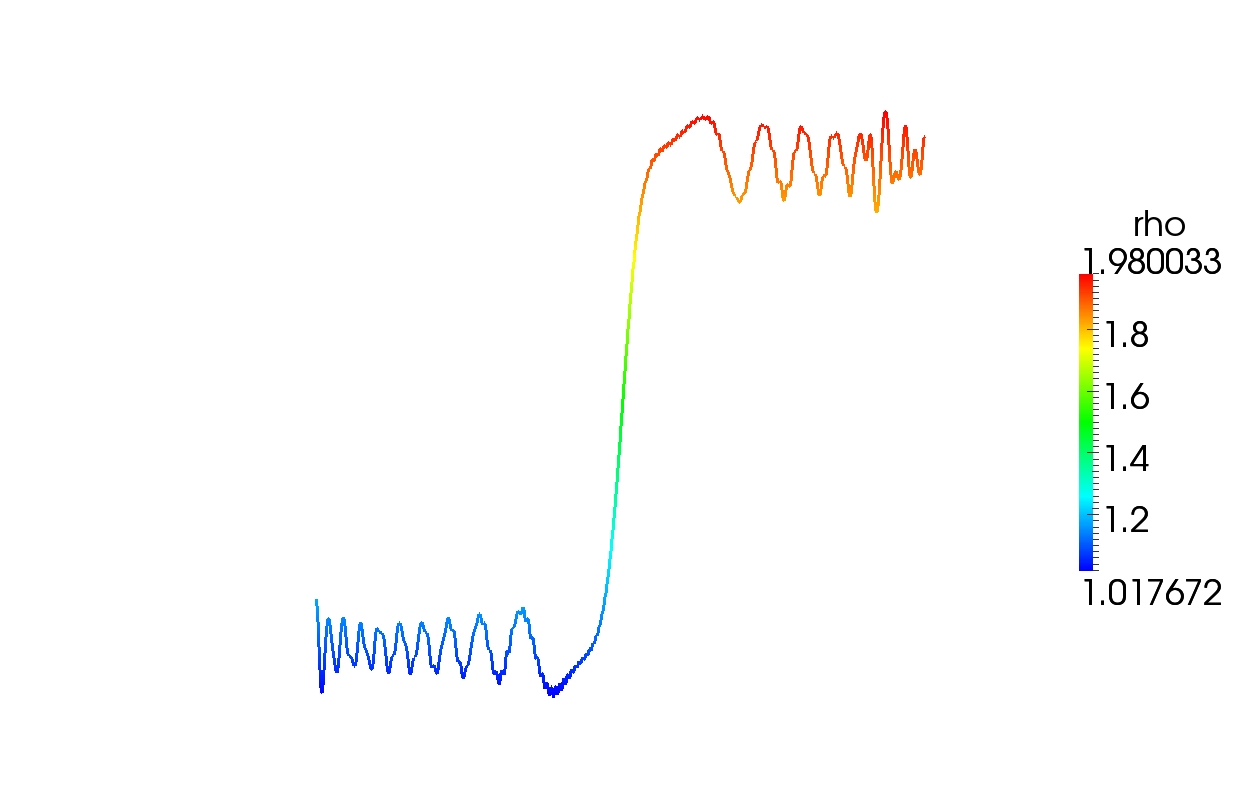}
            } \subfigure[][$t = 0.5$]{ \includegraphics[scale=\figscale,
                width=0.47\figwidth] {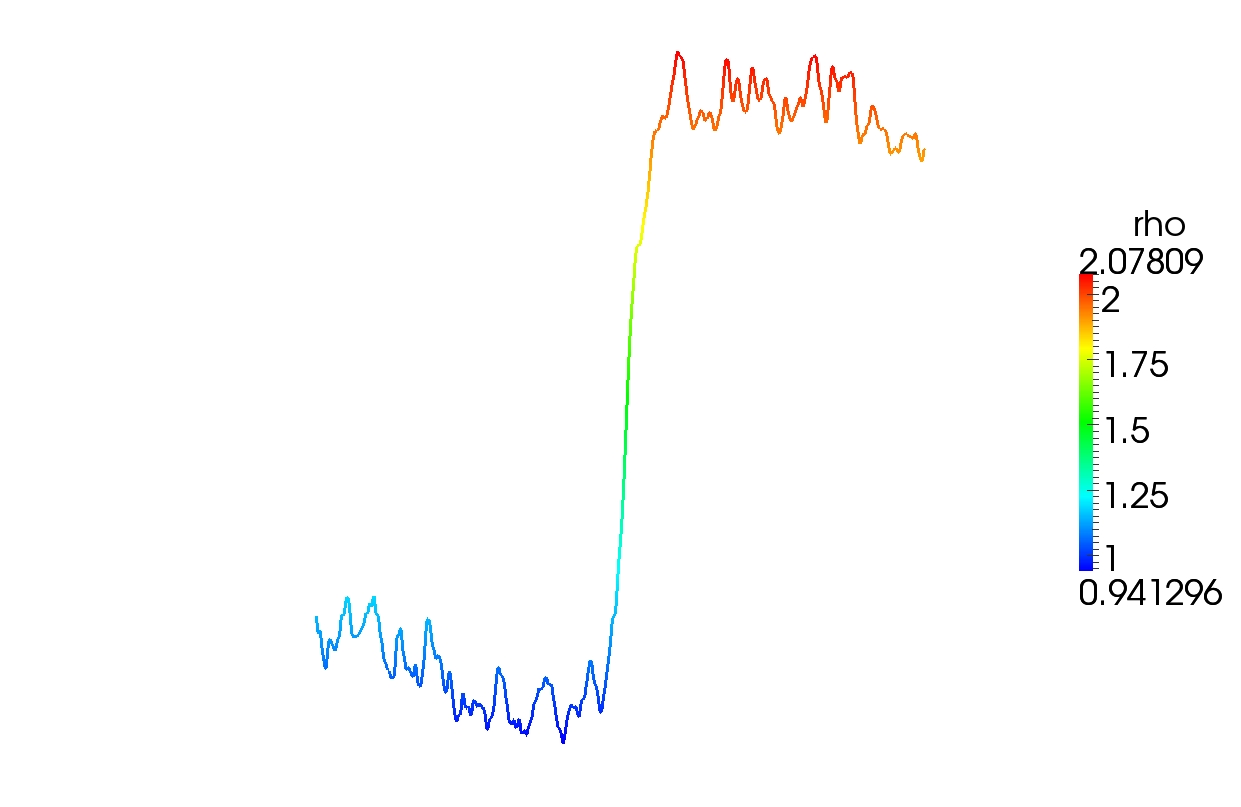}
            } \subfigure[][Conservativity plot]{ \includegraphics[scale=\figscale,
                width=0.47\figwidth]
              {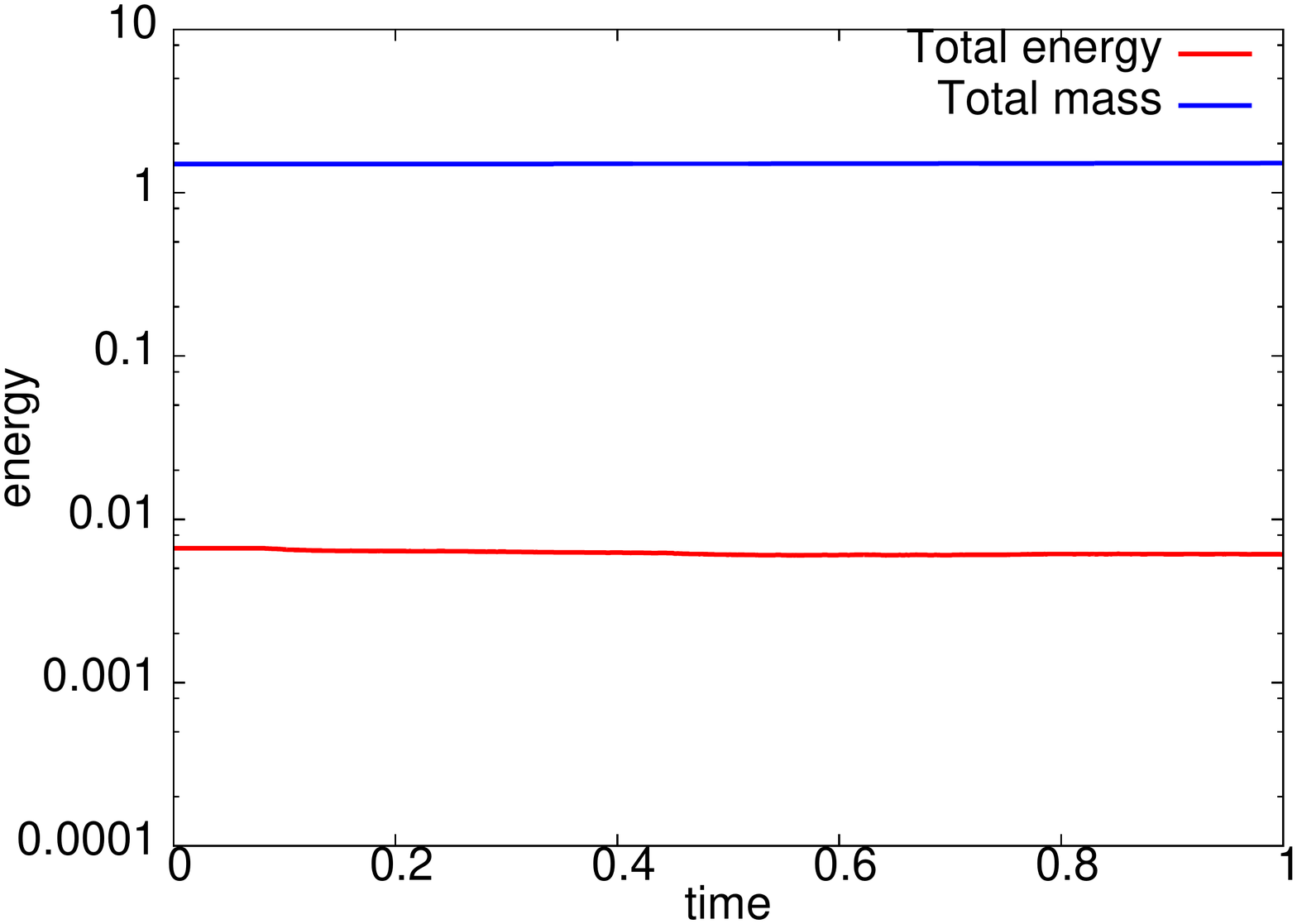} }
          \end{center}
\end{figure}

We take $\gamma = 10^{-4}$, $h = 10^{-4}$ and $k_n = k = 10^{-3}$
for each $n$. Figure \ref{fig:1d-conservative-energy-plot} shows the
energy and mass conservativity of the simulation.

\subsection{Test 2 -- monotone energy dissipation for the \NSK system}
\label{sec:energy-NSK}

In this case we take $\mu > 0$ and study the dissipation property
for the full \NSK system given by (\ref{eq:fully_discrete}). We take
$\W = [0,1]$ and consider the initial conditions
(\ref{eq:step-ics}). We fix $\gamma = 10^{-4}, h = 10^{-4}$ and $k_n =
k = 10^{-3}$. 

\begin{figure}[ht]
  \caption[]
          {
            \label{fig:1d-dissipative-energy-plot-mu-10-7} 
            \ref{sec:energy-NSK} Test 2 -- Numerical experiment
            showing the conservation of mass and dissipation of energy
            for the numerical method proposed in Corollary
            \ref{cor:consistency-and-conservation} for the \NSK
            system. In this test we take $\mu = 10^{-7}$. Notice that
            the energy dissipation allows the \NSK simulation to
            achieve a steady state. At $t=50$ the maximal value of the
            velocity is of magnitude $10^{-5}$. Notice also that $\mu$
            is chosen sufficiently small such that the dynamics are
            comparible with that of Figure
            \ref{fig:1d-conservative-energy-plot} albeit with smeared out oscillations. }
          \begin{center}
            \subfigure[][Initial condition, $t = 0$]{
              \includegraphics[scale=\figscale, width=0.47\figwidth]
                              {./figures/EK-init-cond-test-1.jpg} }
            \subfigure[][$t = 0.01$]{ \includegraphics[scale=\figscale,
                width=0.47\figwidth] {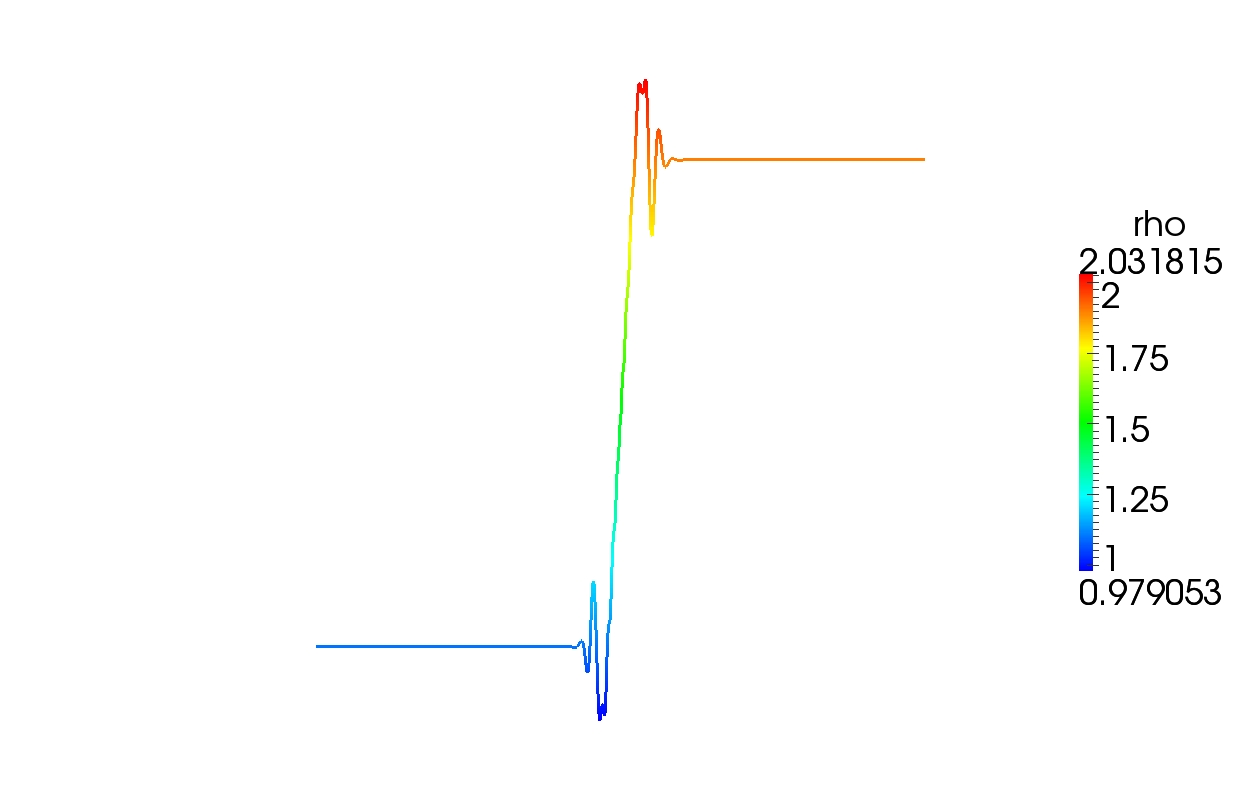}
            } \subfigure[][$t = 0.05$]{ \includegraphics[scale=\figscale,
                width=0.47\figwidth] {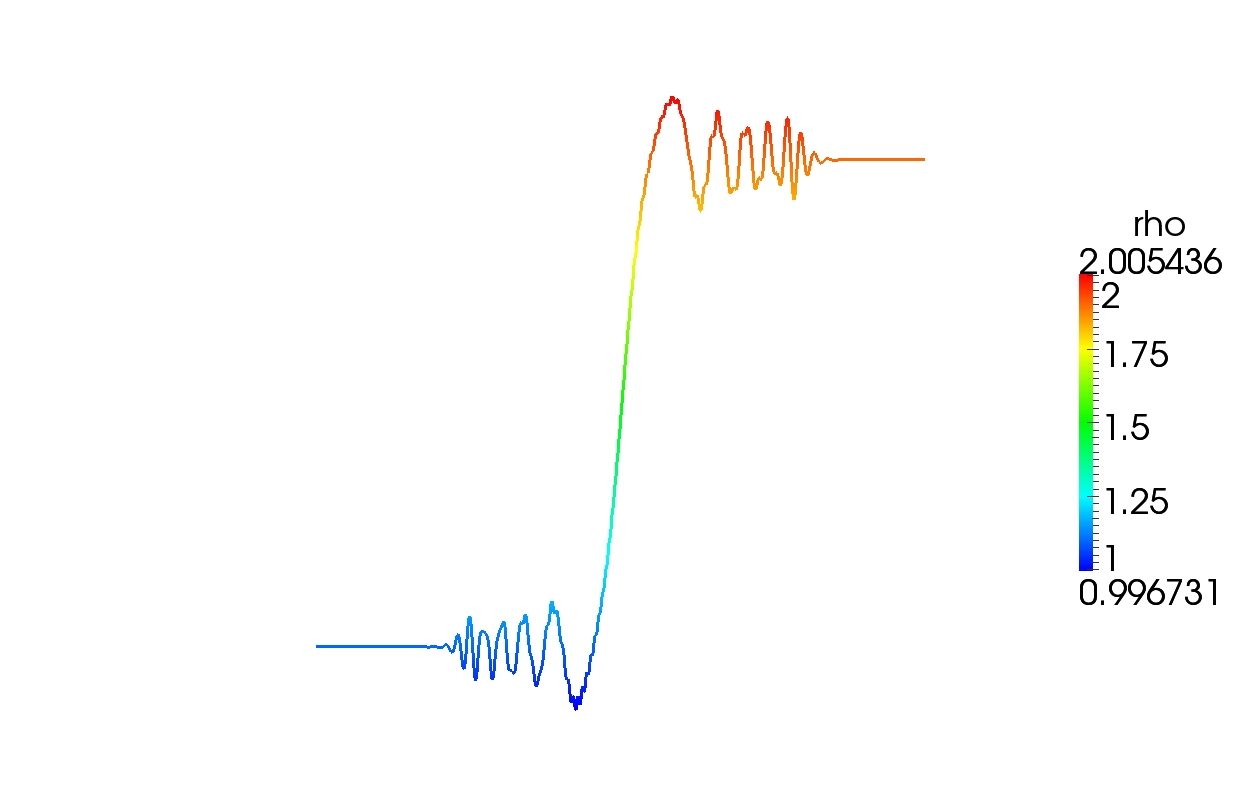}
            } \subfigure[][$t = 0.1$]{ \includegraphics[scale=\figscale,
                width=0.47\figwidth] {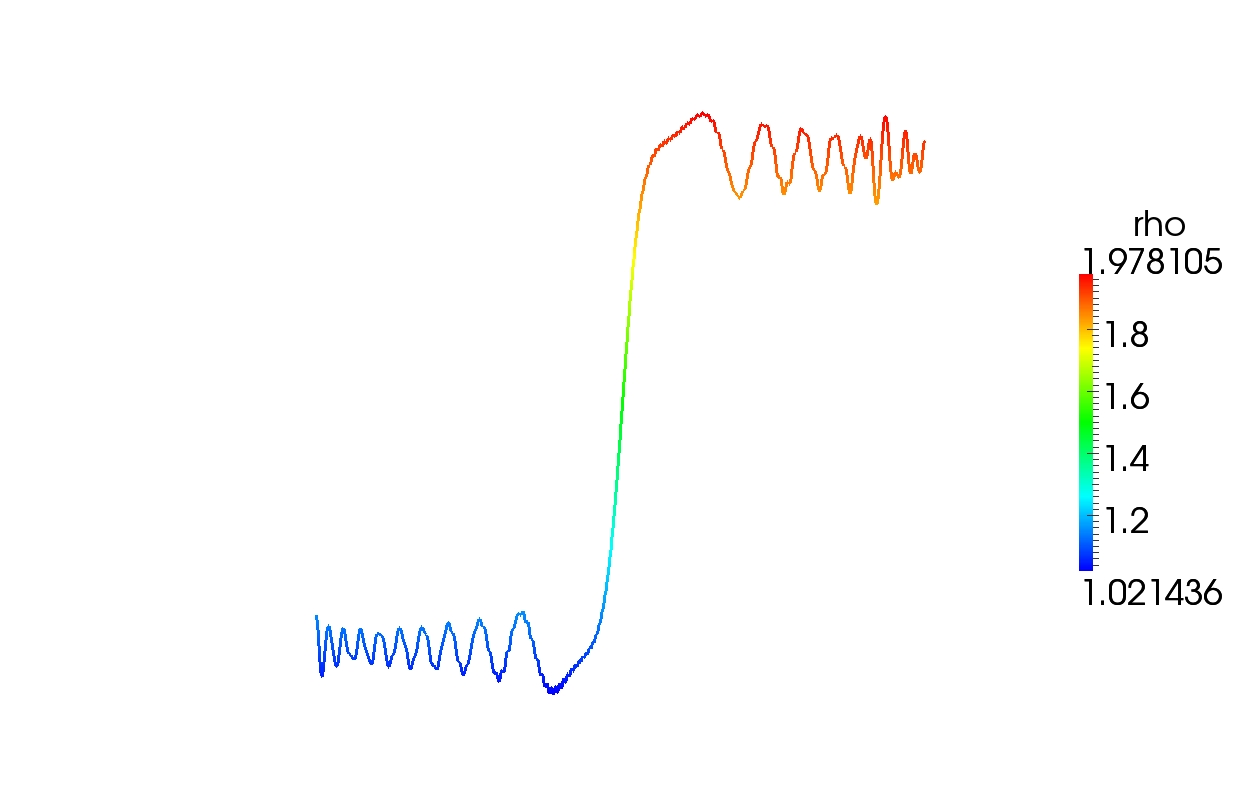}
            } \subfigure[][$t = 0.5$]{ \includegraphics[scale=\figscale,
                width=0.47\figwidth] {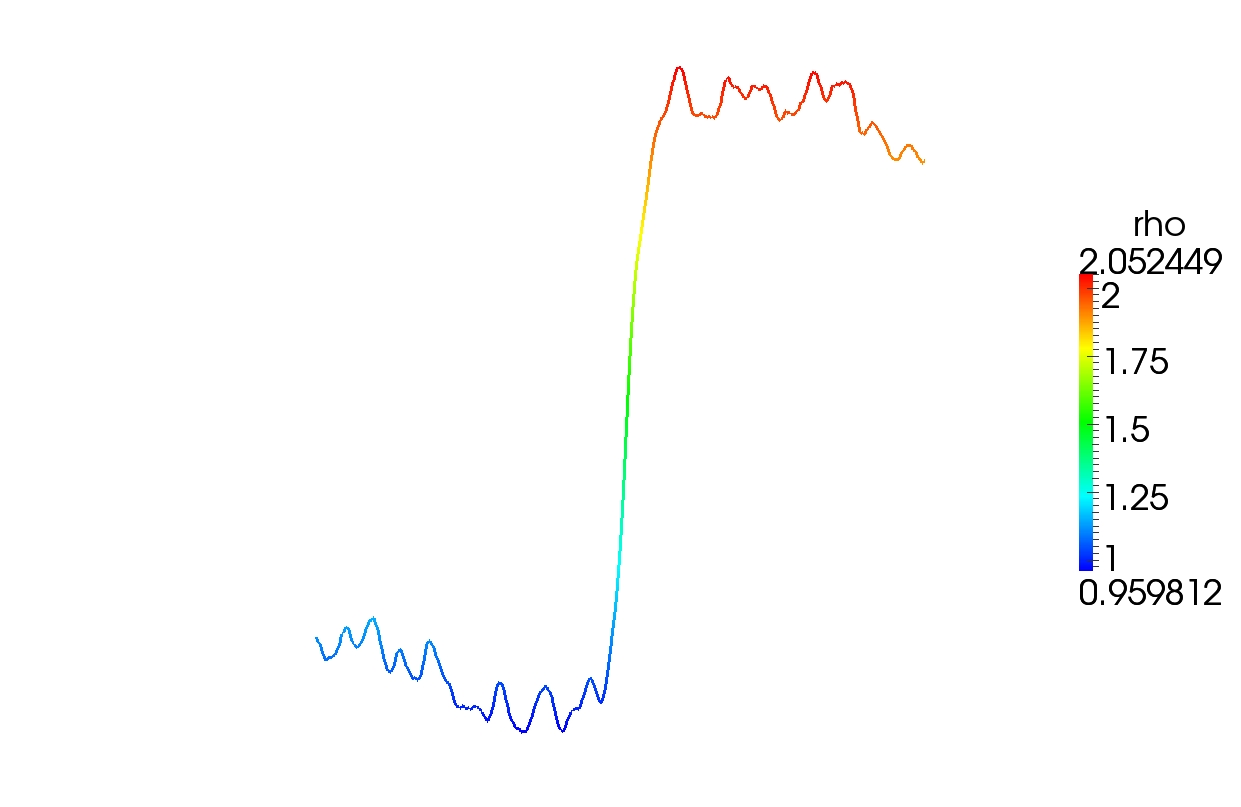}
            } \subfigure[][Conservativity plot]{ \includegraphics[scale=\figscale,
                width=0.47\figwidth]
              {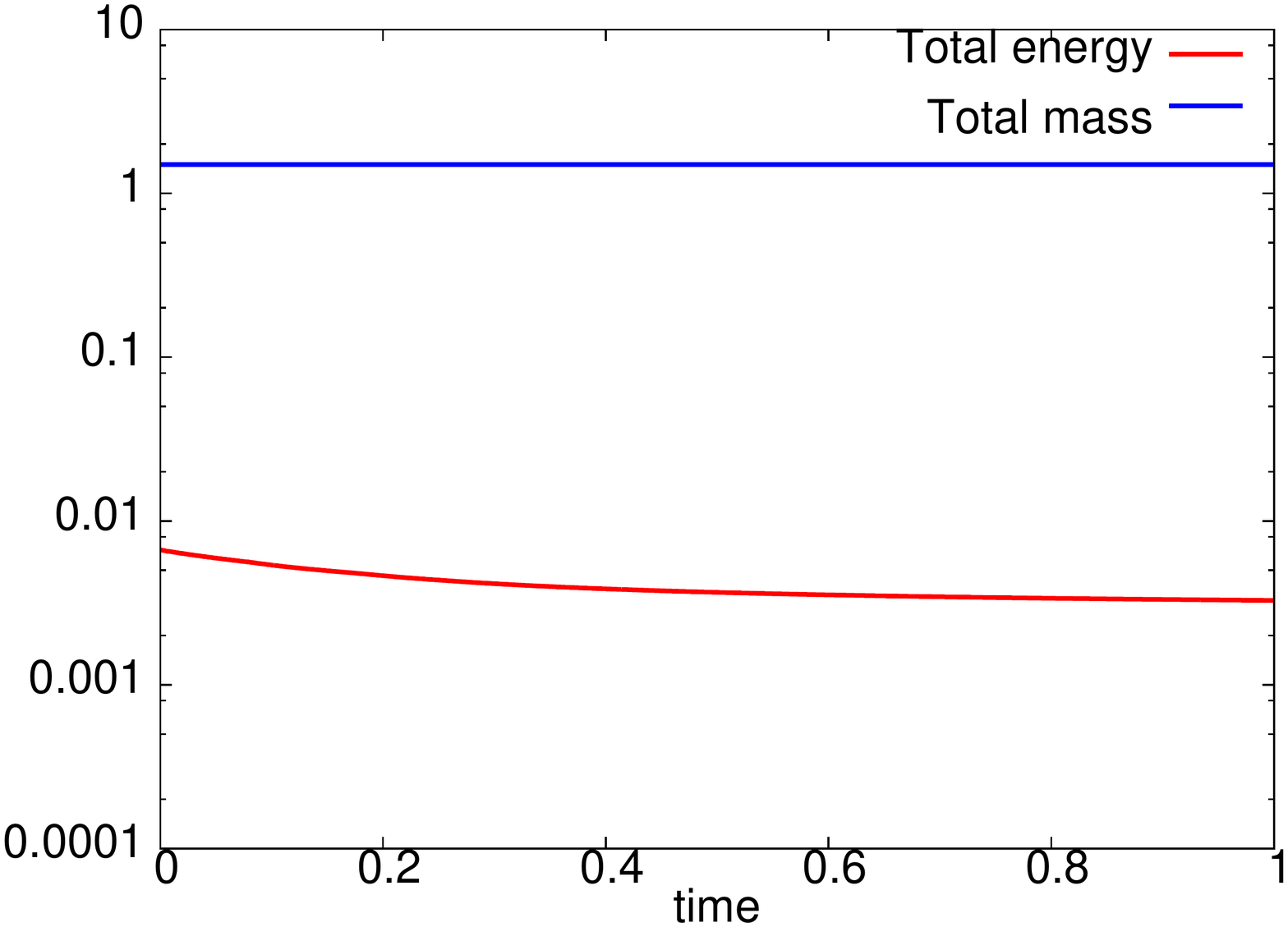} }
          \end{center}
\end{figure}

\begin{figure}[ht]
  \caption[]
          {
            \label{fig:1d-dissipative-energy-plot-mu-10-6} 
            \ref{sec:energy-NSK} Test 2 -- Numerical experiment
            showing the effect of the ratio of viscocity to
            capillarity on the dynamics of the simulation. The
            simulation is the same as in Figure
            \ref{fig:1d-dissipative-energy-plot-mu-10-7} with the
            exception that $\mu = 10^{-6}$. Notice the oscillations
            have become smeared out. The maximal value of velocity is
            of magnitude $10^{-5}$ at $t = 14$.}
          \begin{center}
            \subfigure[][Initial condition, $t = 0$]{
              \includegraphics[scale=\figscale, width=0.47\figwidth]
                              {./figures/EK-init-cond-test-1.jpg} }
            \subfigure[][$t = 0.01$]{ \includegraphics[scale=\figscale,
                width=0.47\figwidth] {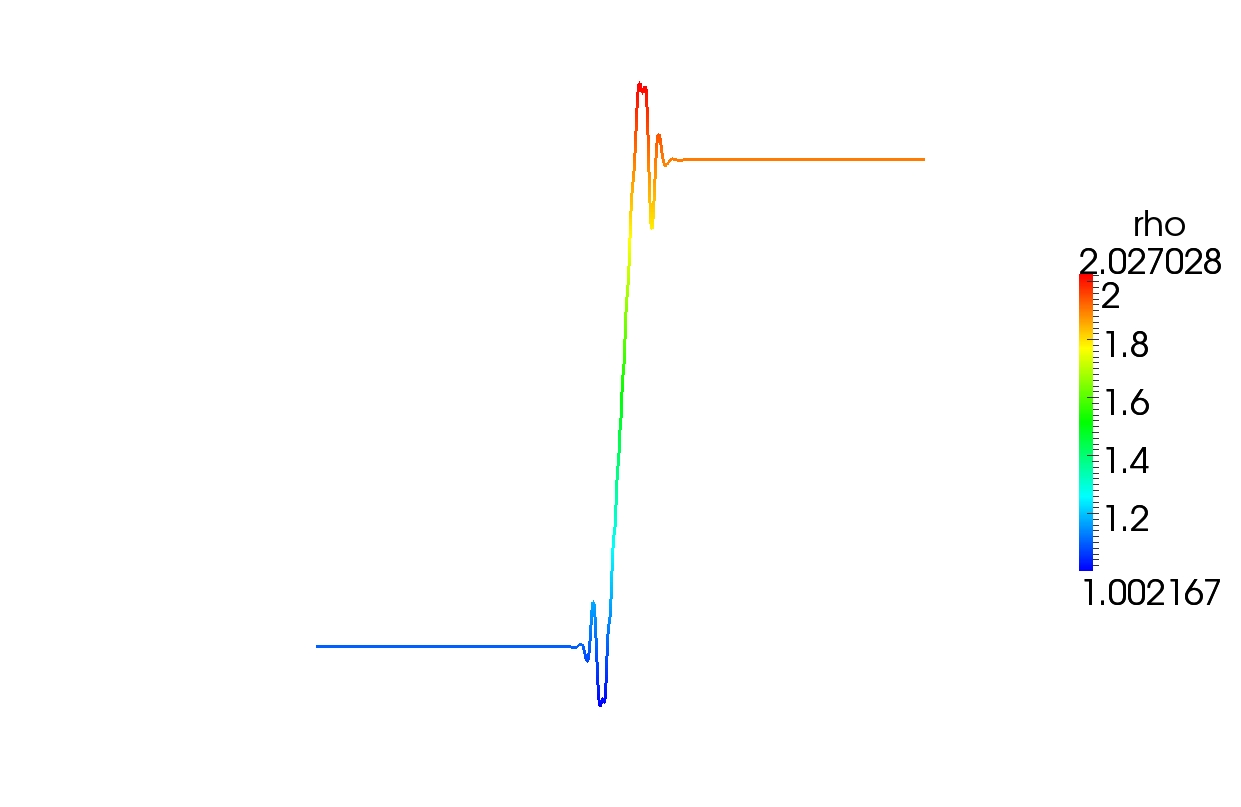}
            } \subfigure[][$t = 0.05$]{ \includegraphics[scale=\figscale,
                width=0.47\figwidth] {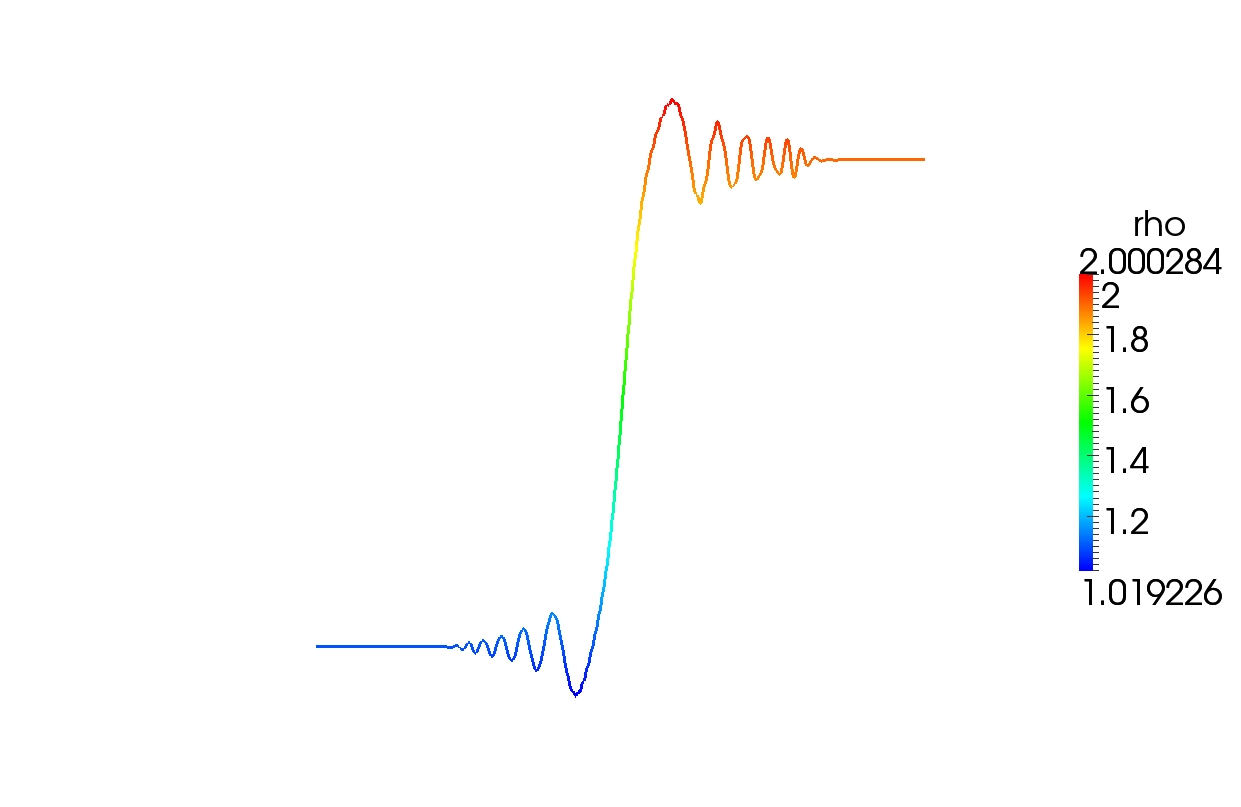}
            } \subfigure[][$t = 0.1$]{ \includegraphics[scale=\figscale,
                width=0.47\figwidth] {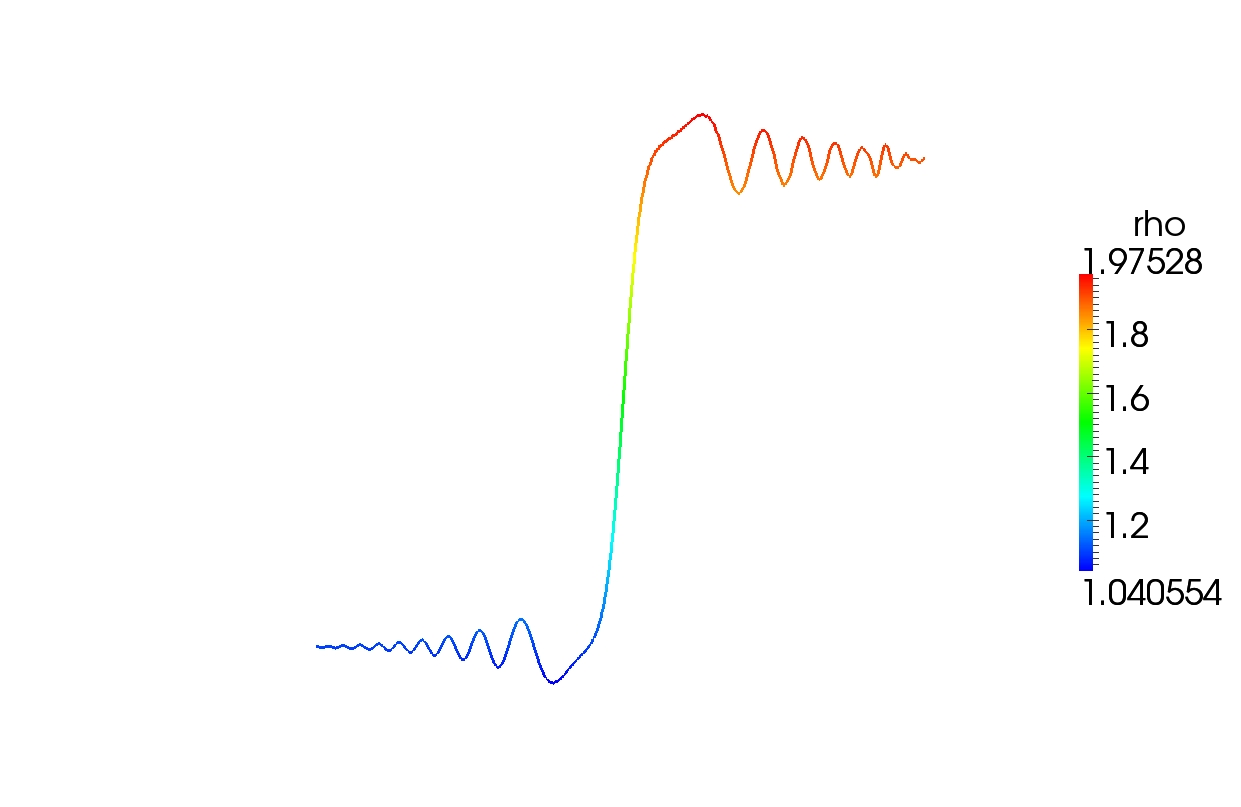}
            } \subfigure[][$t = 0.5$]{ \includegraphics[scale=\figscale,
                width=0.47\figwidth] {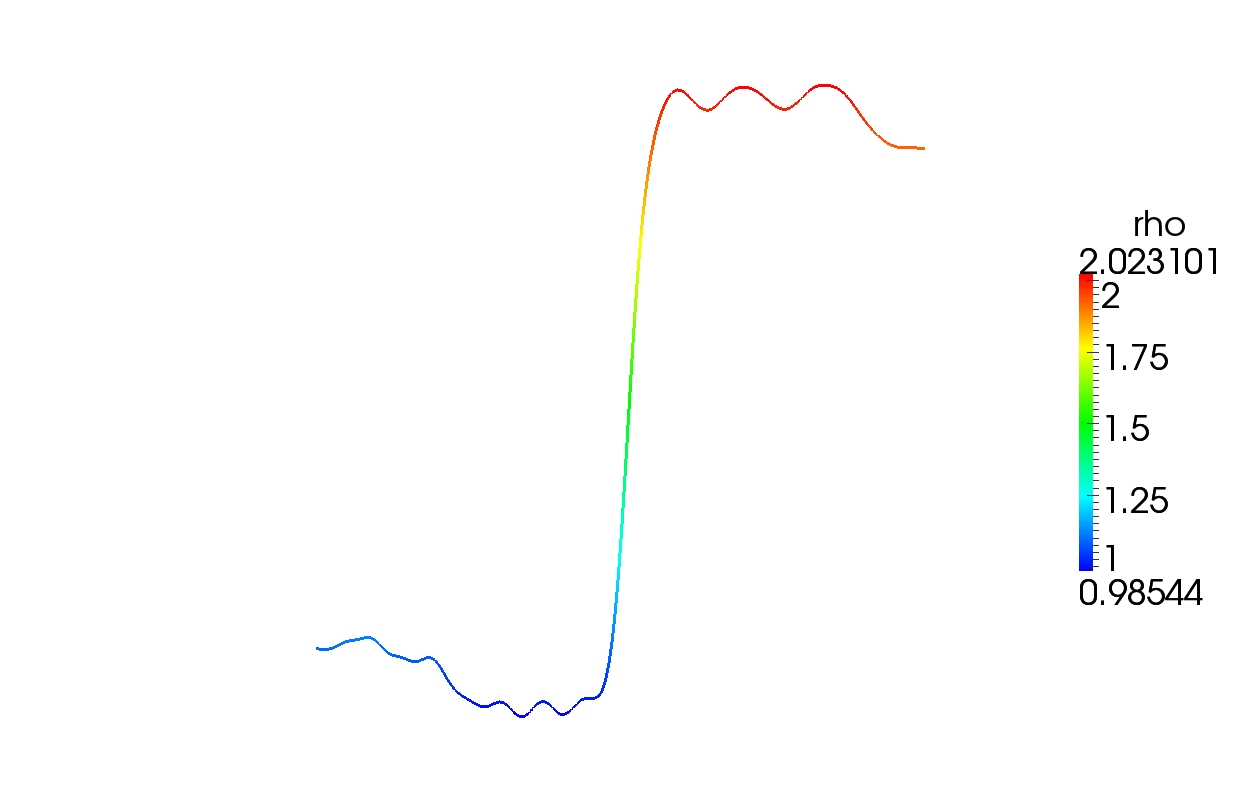}
            } \subfigure[][Conservativity plot]{ \includegraphics[scale=\figscale,
                width=0.47\figwidth]
              {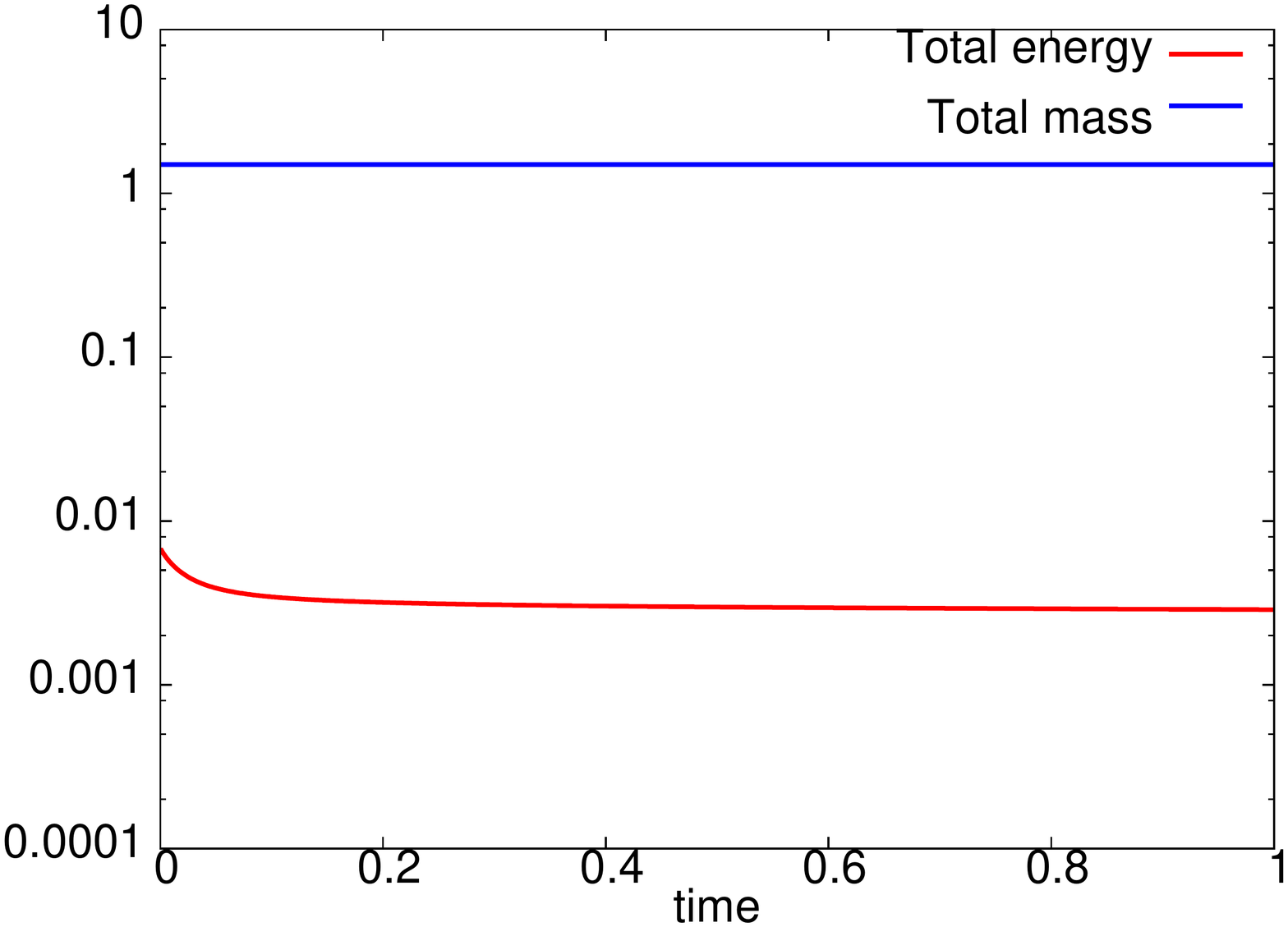} }
          \end{center}
\end{figure}

\begin{figure}[ht]
  \caption[]
          {
            \label{fig:1d-dissipative-energy-plot-mu-10-5} 
            \ref{sec:energy-NSK} Test 2 -- Numerical experiment
            showing the effect of the ratio of viscocity to
            capillarity on the dynamics of the simulation. The
            simulation is the same as in Figure
            \ref{fig:1d-dissipative-energy-plot-mu-10-7} with the
            exception that $\mu = 10^{-5}$. Notice the oscillations
            have become heavily reduced in very short time due to the
            massive dissipation in energy initially. The maximal value of velocity is
            of magnitude $10^{-5}$ at $t = 2$.}
          \begin{center}
            \subfigure[][Initial condition, $t = 0$]{
              \includegraphics[scale=\figscale, width=0.47\figwidth]
                              {./figures/EK-init-cond-test-1.jpg} }
            \subfigure[][$t = 0.01$]{ \includegraphics[scale=\figscale,
                width=0.47\figwidth] {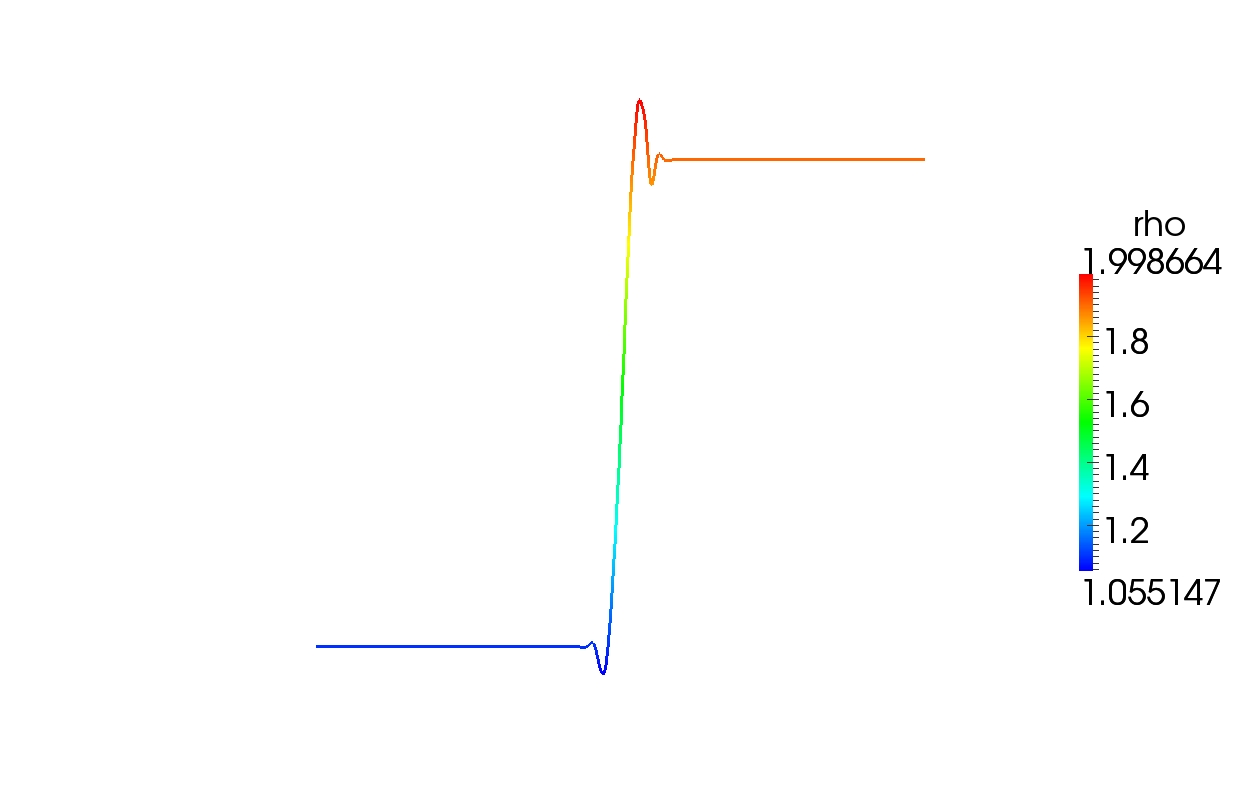}
            } \subfigure[][$t = 0.05$]{ \includegraphics[scale=\figscale,
                width=0.47\figwidth] {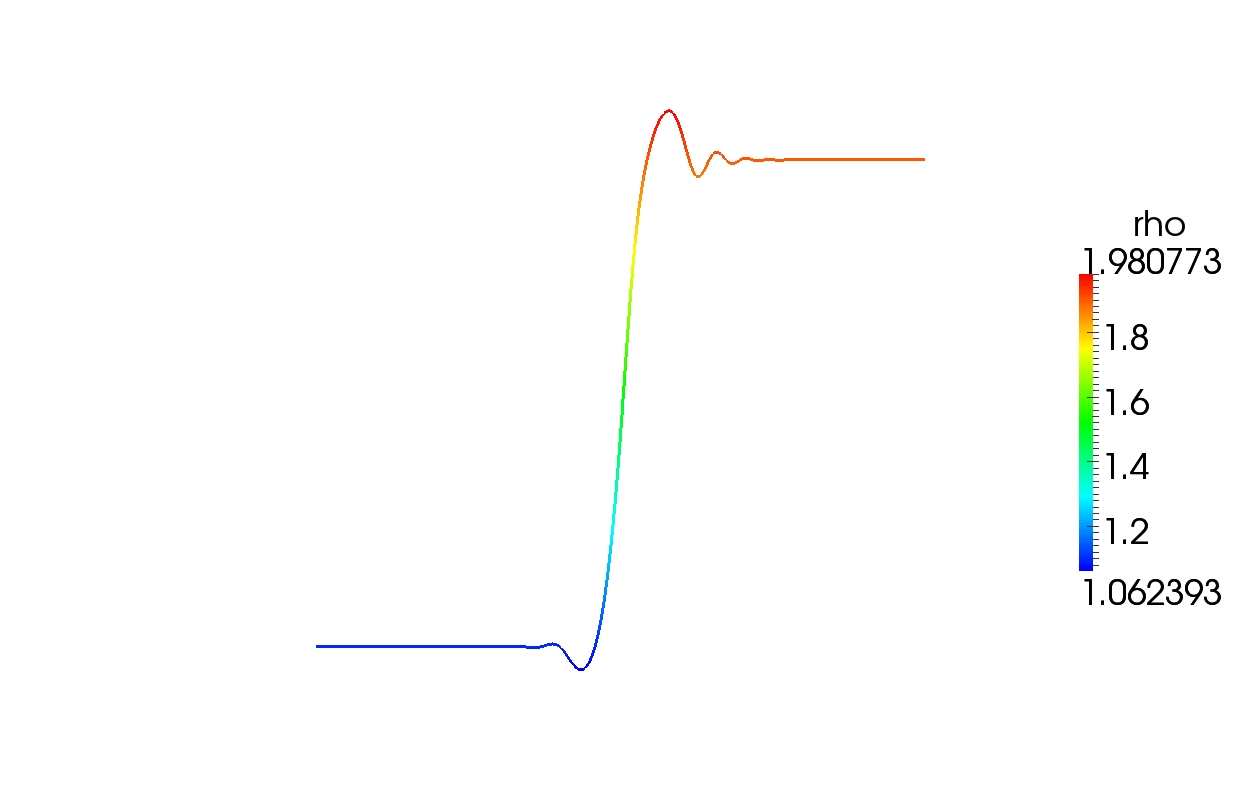}
            } \subfigure[][$t = 0.1$]{ \includegraphics[scale=\figscale,
                width=0.47\figwidth] {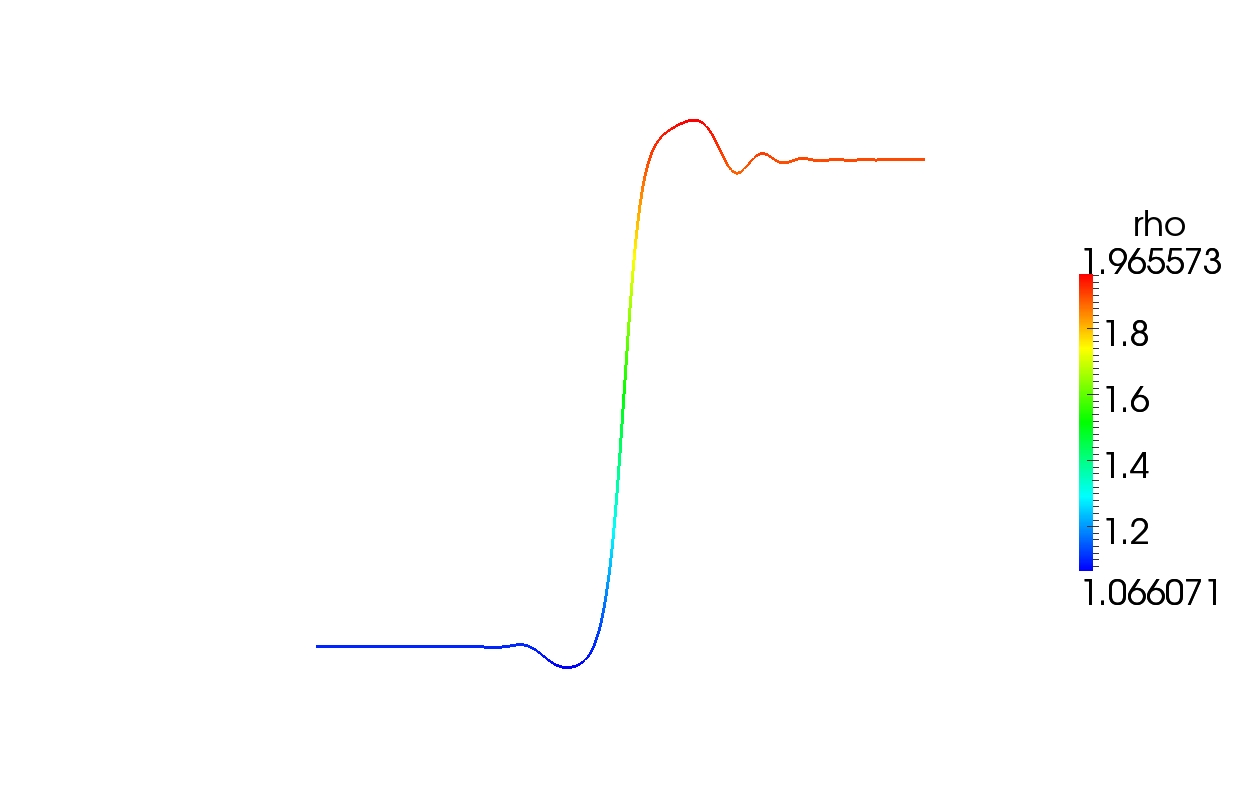}
            } \subfigure[][$t = 0.5$]{ \includegraphics[scale=\figscale,
                width=0.47\figwidth] {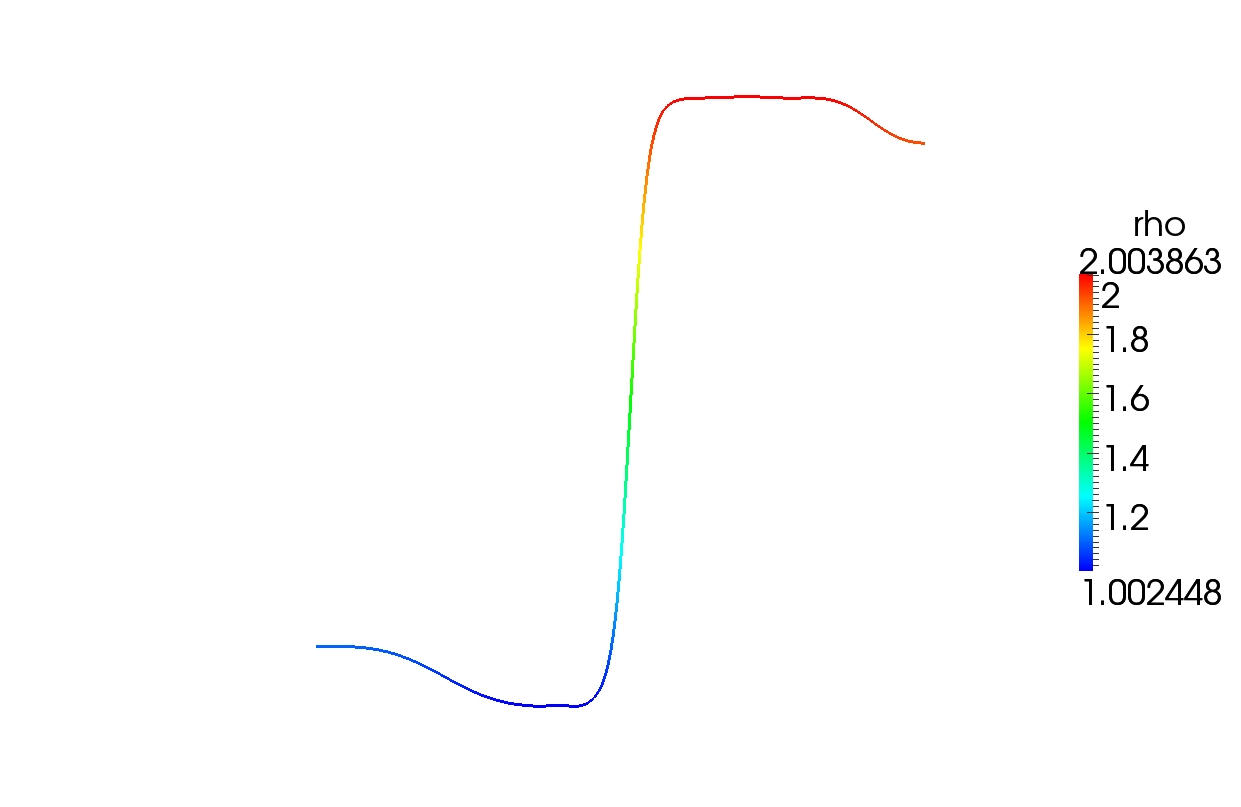}
            } \subfigure[][Conservativity plot]{ \includegraphics[scale=\figscale,
                width=0.47\figwidth]
              {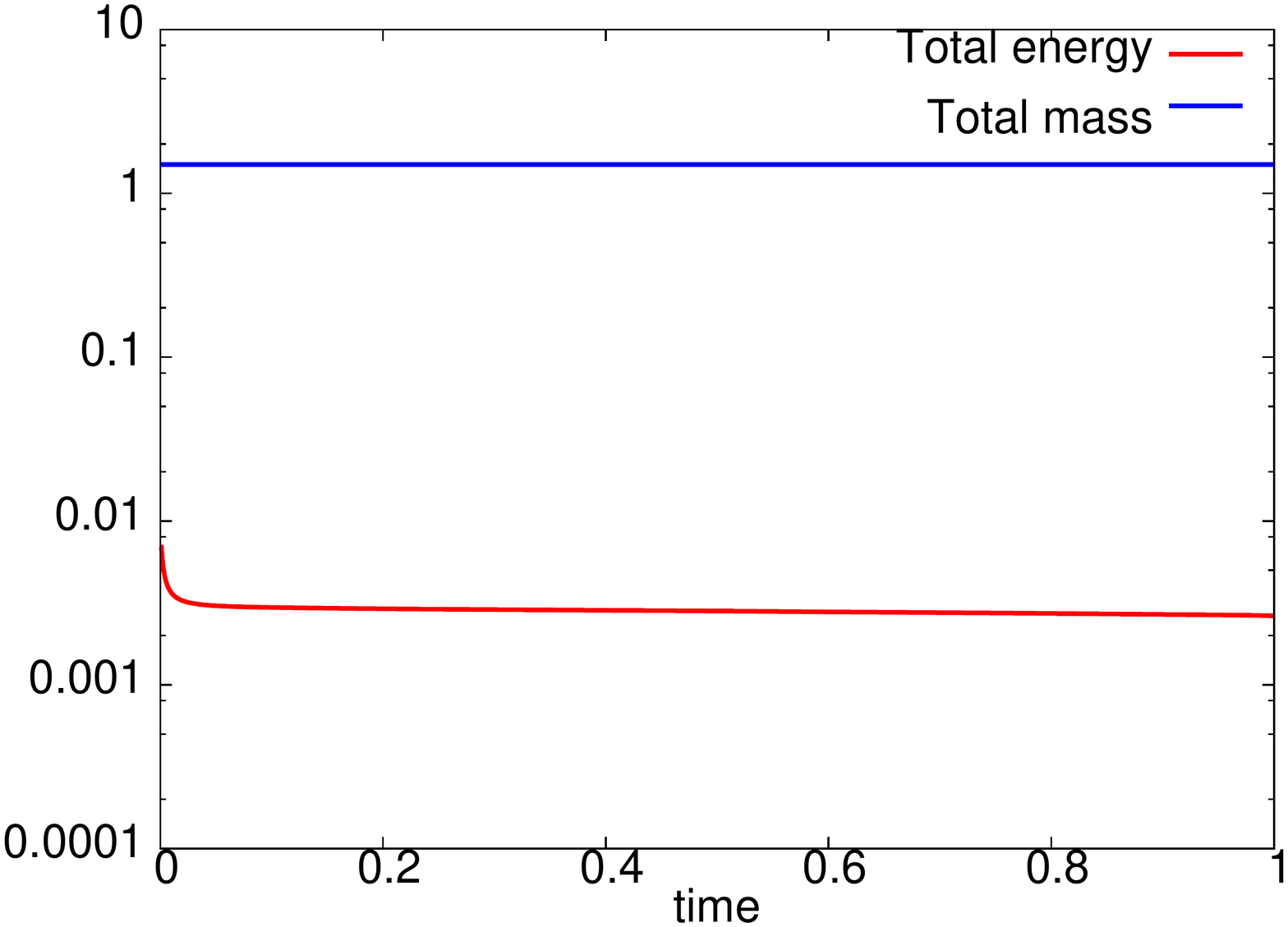} }
          \end{center}
\end{figure}

We test the effect of the ratio of viscocity to capillarity, \ie
$\mu/\gamma$, on the dynamics of the simulation. To that end we run
the simulation for $\mu = 10^{-7}$ (Figure
\ref{fig:1d-dissipative-energy-plot-mu-10-7}), $\mu = 10^{-6}$ (Figure
\ref{fig:1d-dissipative-energy-plot-mu-10-6}) and $\mu = 10^{-5}$
(Figure \ref{fig:1d-dissipative-energy-plot-mu-10-5}).

\subsection{Test 3 -- benchmarking}
\label{sec:benchmark}

In this test we look to benchmark the numerical algorithm against a
steady state solution of the \EK system on the domain $\W = [-1,1]$.

For the double well given by (\ref{eq:double-well}) a steady state
solution to the \EK system is given by 
\begin{gather}
  \rho(x,t) 
  =
  \frac{3}{2} 
  -
  \frac{1}{2}
  \tanh\qp{\frac{x}{2\sqrt{2\gamma}}} 
  \\
  v(x,t) \equiv 0 \Foreach t
\end{gather}
with appropriate initial data. Note that on the boundary $\nabla \rho$
is not zero but of negligable value (for small values of
$\gamma$). Tables \ref{table:p1-gamma-10-4}--\ref{table:p1-gamma-10-6}
detail three experiments aimed at testing the convergence properties
for the scheme for $\gamma = 10^{-4}$ (Table
\ref{table:p1-gamma-10-4}), $\gamma = 10^{-5}$ (Table
\ref{table:p1-gamma-10-5}) and $\gamma = 10^{-6}$ (Table
\ref{table:p1-gamma-10-6}).

\begin{table}[ht]
  \caption{\label{table:p1-gamma-10-4} In this test we benchmark a
    stationary solution of the \EK system using the discretisation
    (\ref{eq:fully_discrete}) with piecewise linear elements ($p =
    1$), choosing $k = 1/N$. This is done by formulating
    (\ref{eq:fully_discrete}) as a system of nonlinear equations, the
    solution to this is then approximated by a Newton method with
    tolerance set at $10^{-15}$. At each Newton step the solution to
    the linear system of equations is approximated using a stabilised
    conjugate gradient iterative solver with an successively
    overrelaxed preconditioner, also set at a tolerance of
    $10^{-15}$. We look at the $\leb{\infty}(0,T; \leb{2}(\W))$ errors
    of the discrete variables $\rho_h$ and $v_h$, and use $e_\rho :=
    \rho - \rho_h$ and $e_v := v - v_h$. In this test we choose
    $\gamma = 10^{-4}$.}
  \begin{center}
    \begin{tabular}{c|c|c|c|c} 
      $N$ & $\Norm{e_\rho}_{\leb{\infty}(\leb{2})}$ & EOC & $\Norm{e_v}_{\leb{\infty}(\leb{2})}$ & EOC \\ 
      \hline  
      32   & 6.258e-3 & 0.000 & 6.194e-4 & 0.000 \\ 
      64   & 3.028e-4 & 4.369 & 4.631e-5 & 3.742 \\ 
      128  & 4.565e-5 & 2.730 & 1.105e-5 & 2.067 \\ 
      256  & 1.155e-5 & 1.983 & 3.691e-6 & 1.582 \\ 
      512  & 2.945e-6 & 1.972 & 9.916e-7 & 1.896 \\ 
      1024 & 7.368e-7 & 2.000 & 2.528e-7 & 1.972 \\
      2048 & 1.842e-7 & 2.000 & 6.324e-8 & 1.999 \\
      4096 & 4.605e-8 & 2.000 & 1.580e-8 & 2.009
  \end{tabular}
  \end{center}
\end{table}

\begin{table}
  \caption{\label{table:p1-gamma-10-5} The test is the same as in
    Table \ref{table:p1-gamma-10-4} with the exception that
    we take $\gamma = 10^{-5}$.}
  \begin{center}
    \begin{tabular}{c|c|c|c|c} 
      $N$ & $\Norm{e_\rho}_{\leb{\infty}(\leb{2})}$ & EOC & $\Norm{e_v}_{\leb{\infty}(\leb{2})}$ & EOC \\ 
      \hline  
      32   & 7.017e-3 & 0.000 & 1.315e-3 & 0.000 \\ 
      64   & 2.469e-3 & 1.506 & 5.819e-4 & 1.176 \\ 
      128  & 4.411e-4 & 2.485 & 7.672e-5 & 2.923 \\ 
      256  & 2.885e-5 & 3.935 & 5.693e-6 & 3.752 \\ 
      512  & 6.5970-6 & 2.129 & 1.295e-6 & 2.136 \\ 
      1024 & 1.668e-6 & 1.984 & 3.228e-7 & 2.004 \\
      2048 & 4.161e-7 & 2.003 & 8.017e-8 & 2.010 \\
      4096 & 1.040e-7 & 2.001 & 2.001e-8 & 2.003
  \end{tabular}
  \end{center}
\end{table}

\begin{table}
  \caption{\label{table:p1-gamma-10-6} The test is the same as in
    Table \ref{table:p1-gamma-10-4} with the exception that
    we take $\gamma = 10^{-6}$.}
  \begin{center}
    \begin{tabular}{c|c|c|c|c} 
      $N$ & $\Norm{e_\rho}_{\leb{\infty}(\leb{2})}$ & EOC & $\Norm{e_v}_{\leb{\infty}(\leb{2})}$ & EOC \\ 
      \hline  
      32   & 1.883e-2 & 0.000 & 1.488e-3 & 0.000 \\ 
      64   & 9.071e-3 & 1.054 & 8.134e-4 & 0.871 \\ 
      128  & 3.807e-3 & 1.253 & 3.820e-4 & 1.090 \\ 
      256  & 1.005e-3 & 1.922 & 9.110e-5 & 2.051 \\ 
      512  & 6.486e-5 & 3.954 & 5.118e-6 & 4.171 \\ 
      1024 & 4.907e-6 & 3.724 & 6.809e-7 & 2.910 \\
      2048 & 1.016e-6 & 2.272 & 1.445e-7 & 2.236 \\
      4096 & 2.439e-7 & 2.059 & 3.446e-8 & 2.068
  \end{tabular}
  \end{center}
\end{table}

\clearpage
\subsection{Test 4 -- simulations for $d=2$ and parasitic currents}
\label{sec:2d-sims}

In this test we consider the case $d=2$. We take $\W = [0,1]^2$ and look at the 
following initial condition 
\begin{equation}
  \label{eq:2d-parabolic-ics}
  \rho_0(\geovec x) 
  =
  \begin{cases}
    2 \text{ if } \qp{x,y} \in [0.3,0.7]^2
    \\
    1 \text{ otherwise }
  \end{cases}
\qquad v_0 \equiv 0
\end{equation}
and examine its evolution.

We expect due to the non-local part of the energy that interfacial
layers of size $\sim \sqrt{\gamma}$ form, see \cite{Ste88,ORS90} for
an energy argument. This process smoothes the profile. 
Moreover, the length of the interface is reduced such that the quadratic ``droplet'' becomes circular.

\begin{figure}[h]
  \caption[]
          {
            \label{fig:2d-solution-plot} 
            Test \ref{sec:2d-sims}. The solution, $\rho_h$ to the \NSK
            system with initial conditions (\ref{eq:2d-parabolic-ics})
            at various values of $t$, overlayed with the velocity
            $\geovec v_h$. Notice that the are no parasitic currents
            appearing in the interfacial layer. The velocity tends to
            zero over the entire domain as time increases. The
            energy--mass plot of the simulation is also given.}
          \begin{center}
            \subfigure[][$t=0.001$]{
              \includegraphics[scale=\figscale, width=0.37\figwidth]
                              {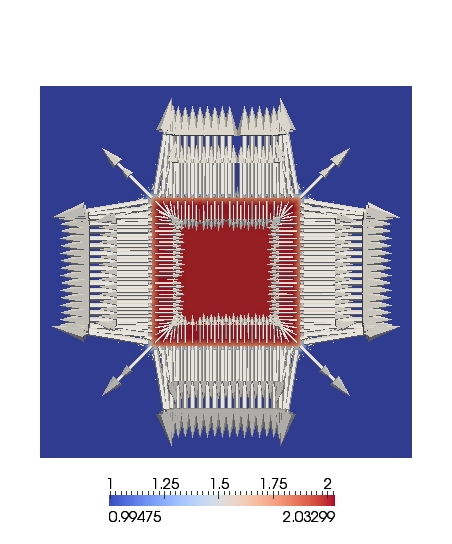}
            }
            \hfill
            \subfigure[][$t=0.1$]{
              \includegraphics[scale=\figscale, width=0.37\figwidth]
                              {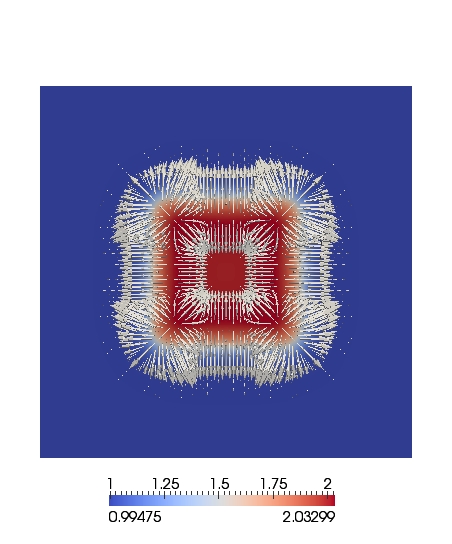}
            }
            \hfill
            \subfigure[][$t=0.25$]{
              \includegraphics[scale=\figscale, width=0.37\figwidth]
                              {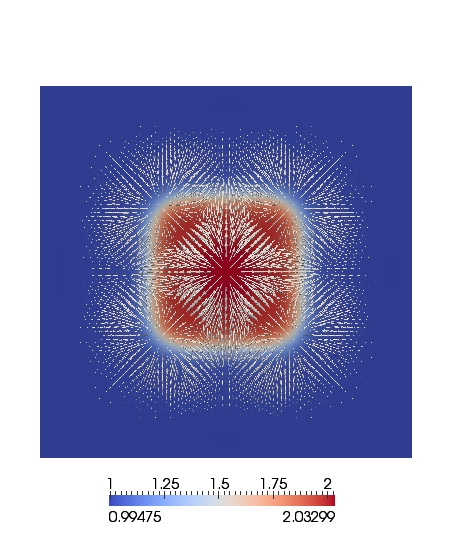}
            }
            \hfill
            \subfigure[][$t=0.5$]{
              \includegraphics[scale=\figscale, width=0.37\figwidth]
                              {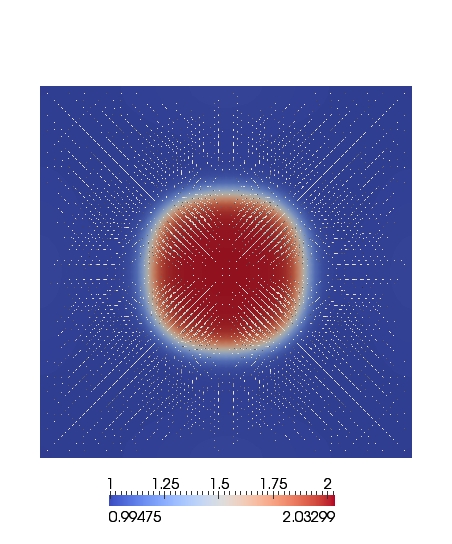}
            }
            \hfill
            \subfigure[][$t=1.4$]{
              \includegraphics[scale=\figscale, width=0.37\figwidth]
                              {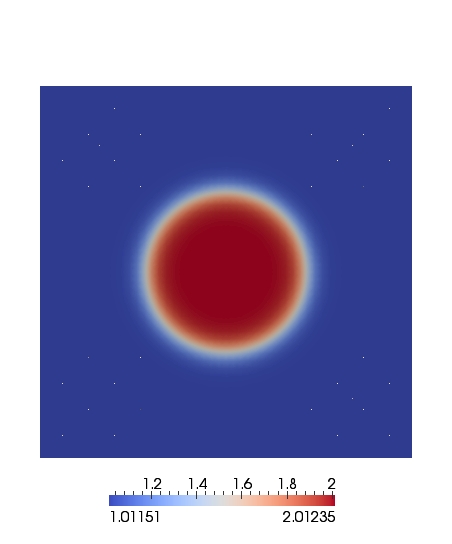}
            }
            \hfill
            \subfigure[][energy--mass]{
              \includegraphics[scale=\figscale, width=0.5\figwidth]
                              {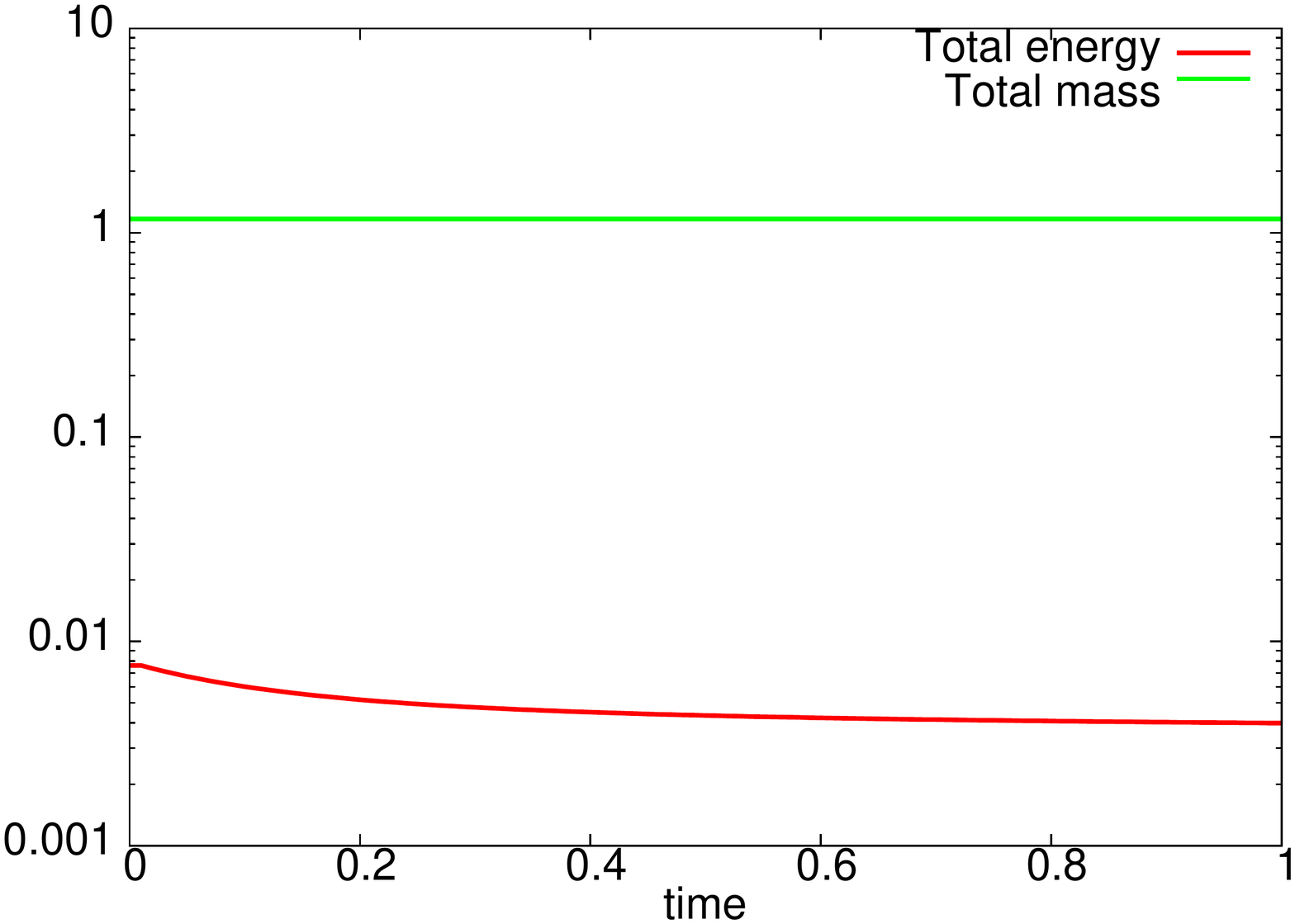}
            }
          \end{center}
\end{figure}

We take $\gamma = \mu = 0.0005$, $h \approx 0.02$ and $k_n = k =
0.001$ for all $n$. Figure \ref{fig:2d-solution-plot} shows the
behaviour of the energy and mass of the numerical solution together
with the solution plot of $\rho_h$ at various times. The solution is
overlayed with the velocity $\geovec v_h$ as a glyph plot.


\clearpage   
\bibliographystyle{alpha}
\bibliography{nskbib,tristansbib}

\end{document}